\documentclass[11pt,a4paper,draft,reqno]{amsart}
\makeatletter
\def\@settitle{\begin{center}\baselineskip14\p@\relax\normalfont\uppercasenonmath\@title\@title\end{center}}
\makeatother
\usepackage{amssymb}
\usepackage{ifthen}
\usepackage{tikz}
\usepackage{tikz-3dplot}
\usepackage{bm}
\usepackage[hyphens]{url} 
\usepackage{amsmath,amssymb,amsfonts,amsthm,amsopn,dsfont,mathrsfs}
\usepackage{esint}
\usepackage{graphicx}
\usepackage{enumitem}
\usepackage{bigints}
\usepackage{fancyhdr}
\usepackage[T1]{fontenc}
\usepackage{mathtools}
\usepackage[utf8]{inputenc}
\usepackage{verbatim} 
\usepackage[backend=bibtex,style=alphabetic,bibencoding=utf8,language=auto,autolang=other,maxnames=10]{biblatex}
\usepackage[a4paper, left=20mm, right=20mm, top=25mm, bottom=25mm, headsep=8mm, footskip=15mm, marginparsep=0mm, marginparwidth=0mm]{geometry}
\theoremstyle{plain}
\newtheorem{thm}{Theorem}[section]

\newtheorem{lem}{Lemma}[section]
\newtheorem{cor}{Corollary}[section]
\theoremstyle{definition}

\newtheorem{dfn}{Definition}[section]
\theoremstyle{remark}
\newtheorem{rem}{Remark}[section]
\numberwithin{equation}{section}

\renewcommand{\r}{\mathbb{R}}
\renewcommand{\P}{\mathcal{P}}
\newcommand{\z}{\mathbb{Z}}
\newcommand{\n}{\mathbb{N}}
\renewcommand{\le}{\leqslant}
\renewcommand{\ge}{\geqslant}

\renewcommand{\H}{\mathcal{H}}
\renewcommand{\L}{\mathcal{L}}
\DeclareRobustCommand{\rchi}{{\mathpalette\irchi\relax}}
\newcommand{\irchi}[2]{\raisebox{\depth}{$#1\chi$}}
\DeclareMathOperator{\card}{card}
\DeclareMathOperator{\Div}{div}
\DeclareMathOperator{\vol}{vol}

\DeclareMathOperator{\dist}{dist}
\DeclareMathOperator{\Lip}{Lip}
\DeclareMathOperator{\supp}{supp}

\DeclareMathOperator{\Int}{Int}

\title{The strong $L^p$-closure of vector fields with finitely many integer singularities on $B^3$}
\author[Riccardo Caniato]{Riccardo Caniato} 
\address{Department of Mathematics \\ ETH Zürich\\ Zürich\\ Switzerland}
\email{riccardo.caniato@math.ethz.ch}
\begin{document}
\begin{abstract}
This paper is aimed to investigate the strong $L^p$-closure $L_{\z}^p(B)$ of the vector fields on the open unit ball $B\subset\r^3$ that are smooth up to finitely many integer point singularities. First, such strong closure is characterized for arbitrary $p\in[1,+\infty)$. Secondly, it is shown what happens if the integrability order $p$ is large enough (namely, if $p\ge 3/2$). Eventually, a decomposition theorem for elements in $L_{\z}^1(B)$ is given, conveying information about the possibility of connecting the singular set of such vector fields by a mass-minimizing, integer $1$-current on $B$ with finite mass. 
\end{abstract}
\maketitle
\tableofcontents
\section{Introduction}
\noindent
\subsection{Useful notation and conventions}
We will always denote by $B$ the open unit ball in $\r^3$ and by $C$ the open unit cube in in $\r^3$. We will indicate the standard euclidean scalar product on $\r^3$ by the symbol "$\cdot$" .\\
Given any $p\in[1,+\infty]$ and $k\in\{0,1,2,3\}$, we denote by $\Omega_{L^p}^k(B)$ the set of $k$-forms 
\[\omega=\sum_{1\le i_1<...<i_k\le 3}\omega_{i_1...i_k}dx^{i_1}\wedge...\wedge d x^{i_k}\]
on $B$ such that $\omega_{i_1...i_k}\in L^p(B)$, for every $1\le i_1<...<i_k\le 3$. The sets $\Omega_{L_{loc}^p}^k(B)$ and is $\Omega_{W^{1,p}}^k(B)$ are defined analogously. We denote by $\Omega^k(B)$ the set of smooth $k$-forms on $B$ and by $\mathcal{D}^k(B)$ the set of compactly supported smooth $k$-forms on $B$. The set of $k$-dimensional currents on $B$ (i.e. the topological dual of $\mathcal{D}^k(B)$) is denoted by $\mathcal{D}_k(B)$.\\
By the symbol "$\sharp$", we indicate the musical isomorphism $\sharp:\Omega_{L_{loc}^1}^1(B)\rightarrow L_{loc}^1(B;\r^3)$ given  by:
\begin{align*}
\omega=\sum_{i=1}^3\omega_idx^i\, \mapsto \, \omega^{\sharp}:=(\omega_1,\omega_2,\omega_3), \qquad \forall\, \omega\in\Omega_{L_{loc}^1}^1(B).
\end{align*}
We denote by "$\flat$" the inverse operator $\flat: L_{loc}^1(B;\r^3)\rightarrow\Omega_{L_{loc}^1}^1(B)$ of the linear isomorphism $\sharp$.\\
Fix any $k\in\{0,1,2,3\}$. For every $k$-tuple $(i_1,..,i_k)$, we let $(\tilde i_1,...,\tilde i_{3-k})$ be the $(3-k)$-tuple such that 
\[dx^{i_1}\wedge...\wedge dx^{i_k}\wedge dx^{\tilde i_1}\wedge...\wedge dx^{\tilde i_{3-k}}=dx^1\wedge dx^2\wedge dx^3.\]
We denote by "$\star$" the Hodge star operator $\star:\Omega_{L_{loc}^1}^k(B)\rightarrow \Omega_{L_{loc}^1}^{3-k}(B)$ given by
\begin{align*}
\omega=\sum_{1\le i_1<...<i_k\le 3}\omega_{i_1...i_k}dx^{i_1}\wedge...\wedge d x^{i_k}\, \mapsto \, \star\,\omega:=\sum_{1\le i_1<...<i_k\le 3}\omega_{i_1...i_k}dx^{\tilde i_1}\wedge...\wedge d x^{\tilde i_{3-k}},
\end{align*} 
for every $\omega\in\Omega_{L_{loc}^1}^k(B)$.\\
Let $\omega\in\Omega_{L^p}^2(B)$. A priori, it is not clear how we could define the "restriction" of $\omega$ to the boundary of an open cube, since no trace operator is available for $L^p$ functions. Still, we can give a meaning to such a restriction for "a.e." cube in $B$, in the following precise sense. \\
Fix any $x_0\in B$. Consider the function $d_{x_0}:B\rightarrow(0,+\infty)$ given by $n_{x_0}(x):=||x-x_0||_{*}$, where, here and throughout this paper, we let
\begin{align*}
||x||_{*}:=\sup_{j=1,...,3}|x_j|, \qquad \mbox{ for every } x\in\r^3.
\end{align*}
Clearly, $n_{x_0}^{-1}(r)=\partial C_r(x_0)=\partial (rC+x_0)$, for every $r\in(0,r_{x_0})$ with $r_{x_0}:=\frac{2}{\sqrt{3}}\dist(x_0,\partial B)$. Moreover, $n_{x_0}$ is a Lipschitz function on $B$, which implies that $dn_{x_0}\in\Omega_{L^{\infty}}^1(B)$. Thus, the form $\omega\wedge dn_{x_0}\in\Omega_{L^p}^3(B)$ is a top dimension $L^p$-form on $B$. Since
\begin{align*}
\int_0^{r_{x_0}}\int_{\partial C_r(x_0)}|\star(\omega\wedge dn_{x_0})|^p\, d\H^2\, dr&=\int_{n_{x_0}^{-1}((0,r_{x_0}))}|\star(\omega\wedge dn_{x_0})|^p\, d\L^3\\
&\le\int_{B}|\star(\omega\wedge dn_{x_0})|^p\, d\L^3<+\infty,
\end{align*}
by Fubini's theorem we conclude that for a.e. $r\in (0,r_{x_0})$ the restriction of the function $\star(\omega\wedge dn_{x_0})$ to $\partial C_r(x_0)$ is a well defined $L^p$-function. For every such $r$, we set
\begin{align*}
i_{\partial C_r(x_0)}^*\omega:=\star(\omega\wedge dn_{x_0})|_{\partial C_r(x_0)}\vol_{\partial C_r(x_0)},
\end{align*}
where $\vol_{\partial C_r(x_0)}$ is the Lipschitz volume form on $\partial C_r(x_0)$. We say that $i_{\partial C_r(x_0)}^*\omega$ is the restriction of $\omega$ to $\partial C_r(x_0)$.\\
The reason behind this notation is the following: given any smooth $2$-form $\omega\in\Omega^2(B)$ and denoting by $i_{\partial C_r(x_0)}:\partial C_r(x_0)\hookrightarrow B$ the canonical inclusion map, the procedure that we have described previously generates exactly the pullback of $\omega$ through the map Lipschitz map $i_{\partial C_r(x_0)}$. 

\noindent
\subsection{Goal of the paper and related literature}
The goal of the present paper is to investigate the set $L_R^p(B)$ of the vector fields in $L^p(B;\r^3)$ which can be strongly $L^p$-approximated by vector fields $X$ that are smooth on $B$ up a finite singular set $\{x_1,...,x_n\}$ and have distributional divergence given by
 \[\Div(X)=\sum_{j=1}^nd_j\delta_{x_j}, \qquad \mbox{ for some } d_1,...,d_n\in\z\]
 (see Definition \ref{definition vector fileds integer valued fluxes}). We will often take advantage of the linear isomorphism between vector fields and $2$-forms in dimension 3 (which is described later, in Remark \ref{R1.3}) in order to reformulate this analytic problem in a more geometric framework and to use the formalism of differential forms to describe it. We aim to do so both because the language of principal abelian bundles is particularly convenient to ease the notation and because one of the main areas of interest for this kind of strong approximation results is the study of the Yang-Mills functional in supercritical dimension (see e.g. \cite{kessel},
 \cite{weak-abelian} and \cite{ym-plateau}).\\
 There is at least another important reason to look at such class of vector fields. Consider the functional spaces
 \[W^{1,2}(B;S^2):=\big\{u\in W^{1,2}(B;\r^3) \mbox{ s.t. } |u(x)|=1, \mbox{ for } \L^3\mbox{-a.e. } x\in B\big\}\]
 and
 \[R^{1,2}(B;S^2):=\big\{u\in W^{1,2}(B;S^2) \mbox{ s.t. } u\in C^{\infty}(B\smallsetminus\{x_1,...,x_n\}) \mbox{ for some } x_1,...,x_n\in B\big\}.\]
 Denoting by $\omega_{S^2}\in\Omega^2(S^2)$ the standard volume form on $S^2$, for every $u\in W^{1,2}(B;S^2)$ we let $D(u)$ be the vector field which is associated to the differential $2$-form $u^*\omega_{S^2}$ through the isomorphism that we have mentioned before. Explicitly, $D(u)$ is given by
 \[D(u)=(u\cdot \partial_{x_2} u\times\partial_{x_3} u,u\cdot \partial_{x_3} u\times\partial_{x_1} u, u\cdot \partial_{x_1} u\times\partial_{x_2} u), \qquad \mbox{ for every } u\in W^{1,2}(B;S^2).\]
The vector field $D(u)$ is called the $D$-field associated to the function $u$. By \cite{bethuel-zheng}, we know that for every $u\in W^{1,2}(B;S^2)$ there exists a sequence $\{u_k\}_{k\in\n}\subset R^{1,2}(B;S^2)$ such that $u_k\rightarrow u$ strongly with respect to the $W^{1,2}$-norm as $k\rightarrow\infty$. Clearly, $D(u_k)\in L_R^1(B)$, for every $k\in\n$. Moreover, it holds that $D(u_k)\rightarrow D(u)$ strongly in $L^1(B;\r^3)$ as $k\rightarrow\infty$.  Hence, it follows that
 \[D(u)\in\overline{L_R^1(B)}^{L^1}, \qquad \mbox{ for every } u\in W^{1,2}(B;S^2).\]
 Thus, studying the structure of the strong $L^1$-closure of vector fields in $L_R^1(B)$ carries important information about the $D$-fields associated to the maps in $W^{1,2}(B;S^2)$, which can be useful tools in the analysis of harmonic maps from $B$ to $S^2$ (see e.g. \cite{brezis-coron-lieb} and \cite{betheul-bresiz-coron}).\\
 For example, we will be able to show that even if the singular set of a map in $W^{1,2}(B;S^2)$ can be arbitrarily wild it's always possible to connect such singularities through a integer $1$-current in $B$ with finite mass (see Theorem \ref{T1.2}). 
 
\subsection{Statement of the main results}
Throughout this paper, we will denote  by $C$ the open unit cube centred at the origin of $\r^3$ and we will let \[C_\varepsilon(x_0):=\varepsilon C+x_0,\]
for every $x_0\in\r^3$ and for every $\varepsilon>0$.
\begin{dfn}
\label{D1.1}
Let $p\in[1,+\infty)$. We denote by $\Omega_{p,\z}^2(B)$ the set of all the $2$-forms $F\in\Omega_{L^p}^2(B)$ such that
\[\frac{1}{2\pi}\int_{\partial C_r(x_0)}i_{\partial C_r(x_0)}^*F\in\z,\]
for every $x_0\in B$ and for $\L^1$-a.e. Lebesgue point of the function $F_{x_0}$ given by
\[F_{x_0}(r):=\frac{1}{2\pi}\int_{\partial C_r(x_0)}i_{\partial C_r(x_0)}^*F, \qquad \mbox{ for every } r\in \bigg(0,\frac{1}{\sqrt{3}}\dist(x_0,\partial B)\bigg).\]
\end{dfn}
\begin{rem}
\label{R1.1}
Let $p\in[1,+\infty)$. We notice that $\Omega_{p,\z}^2(B)$ is strongly $L^p$-closed.\\
Indeed, pick any $F\in\Omega_{p,\z}^2(B)$ and consider any sequence $\{F_k\}_{k\in\n}\subset\Omega_{p,\z}^2(B)$ such that $F_k\xrightarrow{k\rightarrow\infty}F$ strongly in $L^p$. Fix any $x_0\in B$, set $r_0:=\dist(x_0,\partial B)$ and define the sets
\[E_k:=\bigg\{r\in\big(0,r_0\big) \mbox{ s.t. } r \mbox{ is a Lebesgue point for } (F_j)_{x_0}\mbox{ and } \frac{1}{2\pi}\int_{\partial C_r(x_0)}i_{\partial C_r(x_0)}^*F_k\notin\z\bigg\},\]
\[E_{\infty}:=\bigg\{r\in\big(0,r_0\big) \mbox{ s.t. } r \mbox{ is a Lebesgue point for } F_{x_0} \mbox{ and } \frac{1}{2\pi}\int_{\partial C_r(x_0)}i_{\partial C_r(x_0)}^*F\notin\z\bigg\},\]
\[E:=\Bigg(\bigcup_{k\in\n}E_k\Bigg)\cup E_{\infty},\]
where $(F_j)_{x_0}$ and $F_{x_0}$ are defined as in Definition \ref{D1.1}.\\
Since by hypothesis $\L^1(E_k)=0$ for every $k\in\n$ and $\L^1(E_{\infty})=0$, it holds that $\L^1(E)=0$. Notice that
\begin{align*}\bigintssss_0^{r_0}\left|\int_{\partial C_r(x_0)}i_{\partial C_r(x_0)}^*F_k-\int_{\partial C_r(x_0)}i_{\partial C_r(x_0)}^*F\right|dr&\le\bigintssss_0^{r_0}\int_{\partial C_r(x_0)}|F_k-F|\, d\H^2dr\\
&=\int_{C_{r_0}(x_0)}|F_k-F|\, d\L^3\le||F_k-F||_{L^1(B)}\xrightarrow{k\rightarrow\infty}0.
\end{align*}
Hence, there exists a subsequence $\big\{F_{k_h}\big\}_{h\in\n}$ of $\{F_k\}_{k\in\n}$ such that
\begin{align}\int_{\partial C_r(x_0)}i_{\partial C_r(x_0)}^*F_{k_h}\xrightarrow{h\rightarrow\infty}\int_{\partial C_r(x_0)}i_{\partial C_r(x_0)}^*F, \qquad \mbox{ for a.e. } r\in(0,r_0).\end{align}
Thus, there exists a subset $G\subset(0,r_0)$ such that $\L^1(G)=0$ and (1.1) holds for every $r\in(0,r_0)\smallsetminus G$. Define $\tilde E:=E\cup G$. By assumption, $\L^1(\tilde E)=0$ and for every $r\in (0,r_0)\smallsetminus\tilde E$ both (1.1) and
\[\int_{\partial C_r(x_0)}i_{\partial C_r(x_0)}^*F_{k_h}\in\z, \qquad \mbox{ for every } k\in\n,\]
hold. Then, for every $r\in\tilde E$ 
\[\int_{\partial C_r(x_0)}i_{\partial C_r(x_0)}^*F\in\z\]
and the above statement follows. 
\end{rem}
\begin{dfn}
Let $p\in[1,+\infty)$. We say that a $2$-form $F\in\Omega_{L^p}^2(B)$ is \textbf{smooth away from $2\pi$-integer point singularities} if
\begin{enumerate}	
\item $F\in C^{\infty}\big(B\smallsetminus\{x_1,...,x_m\}\big)$, for some $m\in\n$ and some set of $m$ points $\{x_1,...,x_m\}\subset B$;
\item \[dF=\Bigg(\sum_{j=1}^m2\pi h_j\delta_{x_j}\Bigg)dx^1\wedge dx^2\wedge dx^3, \qquad \mbox{ in } \mathcal{D}'(B),\]	for some integer values $h_1,...,h_m\in\z$.
\end{enumerate}
We denote by $\Omega_{p,R}^2(B)$ the subset of such $2$-forms in $\Omega_{L^p}^2(B)$. 
\end{dfn}
\begin{rem}
\label{R1.2}
Let $p\in[1,+\infty)$. We notice that $\Omega_{p,R}^2(B)\subset\Omega_{p,\z}^2(B)$. In order to see that, pick any $F\in\Omega_{p,R}^2(B)$, $x_0\in B$ and a Lebesgue point $r$ of the function $F_{x_0}$ (see Definition \ref{D1.1}). We notice that
\[\frac{1}{2\pi}\int_{\partial C_r(x_0)}i_{\partial C_r(x_0)}^*F=\frac{1}{2\pi}\int_{\partial C_r(x_0)}\big(X\cdot\nu_{\partial C_r(x_0)}\big)\, d\H^2,\]
where $X$ is the $L^p$-vector field on $B$ given by $X:=(\star F)^{\sharp}$. We consider the family of Lipschitz and compactly supported functions on $B$ given by
\[\varphi_{\varepsilon}(x):=\begin{cases}1 & \mbox{ for every } x\in C_{r-\varepsilon}(x_0),\\ \displaystyle{\frac{r+\varepsilon}{2\varepsilon}-\frac{|x-x_0|}{2\varepsilon}} & \mbox{ for every } x\in C_{r+\varepsilon}(x_0)\smallsetminus\overline{C_{r-\varepsilon}(x_0)},\\0 & \mbox{ for every } x\in B\smallsetminus\overline{C_{r+\varepsilon}(x_0)}.\end{cases}\]
By direct computation, we get that
\[\nabla\varphi_{\varepsilon}(x)=\begin{cases}0 & \mbox{ for every } x\in C_{r-\varepsilon}(x_0),\\ \displaystyle{-\frac{1}{2\varepsilon}\frac{x-x_0}{|x-x_0|}} & \mbox{ for every } x\in C_{r+\varepsilon}(x_0)\smallsetminus\overline{C_{r-\varepsilon}(x_0)},\\0 & \mbox{ for every } x\in B\smallsetminus\overline{C_{r+\varepsilon}(x_0)}.\end{cases}\]
Hence, since $r$ is a Lebesgue point of $F_{x_0}$, we get that
\begin{align*}\frac{1}{2\pi}\left<\Div(X),\varphi_{\varepsilon}\right>&=-\frac{1}{2\pi}\left<X,\nabla\varphi_k\right>=-\frac{1}{2\pi}\int_{B}(X\cdot\nabla\varphi_k)\, d\L^3\\
&=\frac{1}{2\pi}\bigg(\frac{1}{2\varepsilon}\int_{r-\varepsilon}^{r+\varepsilon}\int_{\partial C_s(x_0)}\big(X\cdot\nu_{\partial C_s(x_0)}\big)\, d\H^2ds\bigg)\\
&=\frac{1}{2\pi}\bigg(\frac{1}{2\varepsilon}\int_{r-\varepsilon}^{r+\varepsilon}F_{x_0}(s)\, ds\bigg)\xrightarrow{\varepsilon\rightarrow 0^+}\frac{1}{2\pi}F_{x_0}(r)=\frac{1}{2\pi}\int_{\partial C_r(x_0)}i_{\partial C_r(x_0)}^*F.
\end{align*}
Moreover, we easily see that 
\[\frac{1}{2\pi}\left<\Div(X),\varphi_{\varepsilon}\right>\xrightarrow{\varepsilon\rightarrow 0^+}\sum_{j\in m_{x_0,r}}h_j\in\z\]
where $m_{x_0,r}\subset m$ is the subset of $m$ given by all the indexes $j$ such that $x_j\in C_r(x_0)$.\\
By uniqueness of the limit, the statement follows.  
\end{rem}
\noindent
The main part of the present paper concerns the proof of the following strong approximation theorem for elements in $\Omega_{p,\z}^2(B)$. Although we point out that this result was previously announced in \cite{kessel}, the proof that we present here is meant to solve some important issues that have never been clarified before. 
\begin{thm}
\label{T1.1}
Let $F\in\Omega_{p,\z}^2(B)$, with $p\in[1,+\infty)$. Then, there exists a sequence $\{F_k\}_{k\in\n}\in\Omega_{p,R}^2(B)$ such that
\[\big|\big|F_k-F\big|\big|_{L^p(B)}\rightarrow 0, \qquad \mbox{ as } k\rightarrow+\infty\]
\end{thm}
\noindent
The characterization of the strong $L^p$-closure of the space $\Omega_{p,\z}^2(B)$ that we mentioned in the previous subsection is an immediate byproduct of Theorem \ref{T1.1}, Remark \ref{R1.1} and Remark \ref{R1.2}. 
\begin{cor}
Let $p\in[1,+\infty)$. Then, 
\[\overline{\Omega_{p,R}^2(B)}^{L^p}=\Omega_{p,\z}^2(B).\]
\end{cor}
\noindent
The proof of Theorem \ref{T1.1} is quite technical and, for the reader's convenience, we sketch here its outline in order to introduce the main ideas on which it is based.\\ 
In section 2, we give some very basic definitions and notions concerning Sobolev principal $G$-bundles. Some kind of acquaintance about these topics is required to understand the vocabulary that will be used intensively in the subsequent sections. A reader who is familiar with this kind of geometric language (in particular, with the theory of Sobolev principal abelian bundles) might want either to skip this section or to read it just to fix the notation.\\
In section 3, we prove the purely technical Lemma \ref{L3.1}. The statement means, roughly speaking, that given any $F\in\Omega_{p,\z}^2(B)$ we can choose a suitable collection $\mathscr{F}=\{C_{\varepsilon}\}_{\varepsilon\in R\subset (0,1)}$ of families of open cubes in $B$ (i.e. any $C_{\varepsilon}$ is a family of open cubes in $B$), each of which will be called a "cubic decomposition" of $B$, such that 
\begin{enumerate}
\item $0$ is a right accumulation point for $R$;
\item every cube in $C_{\varepsilon}$ has side-length $\varepsilon$ and $\L^3(B\smallsetminus C_{\varepsilon})\rightarrow 0$ as $\varepsilon\rightarrow 0^+$ in $R$;
\item the restriction of $F$ to the boundary of every cube contained in each family $C_{\varepsilon}$ belonging to the collection $\mathscr{F}$ has integral in $2\pi\z$;
\item if $\partial C_{\varepsilon}$ is the union of all the boundaries of the cubes in $C_{\varepsilon}$, then the distance of $i_{\partial C_{\varepsilon}}^*F$ in $L^p(\partial C_{\varepsilon})$ from the piecewise constant $2$-form given by taking (in a suitable way that will be specified later) the integral mean of $F$ on each cube contained in $C_{\varepsilon}$ is $\displaystyle{o(\varepsilon^{-1/p})}$ as $\varepsilon\rightarrow 0^+$ in $R$. 
\end{enumerate}
By this point, we can introduce the main idea of the proof. A là Calderón-Zygmund, for every $\varepsilon\in E$ we split the cubes composing $C_{\varepsilon}$ in:
\begin{enumerate}
\item good cubes, whenever the integral of the restriction of $F$ to the boundary of such cubes is less than $2\pi$ and thus, by Lemma \ref{L3.1}, it is actually $0$;
\item bad cubes, otherwise. 
\end{enumerate}
The last Lemma \ref{L3.2} clarifies the first of many advantages that come from having chosen the cubic decompositions so carefully: the volume of the union of the bad cubes vanishes at the limit as $\varepsilon\rightarrow 0^+$ in $E$.\\
The whole Section 4 is aimed to show that, for every fixed $\varepsilon>0$, we can build a smooth connection on the boundary of each cube composing $C_{\varepsilon}$ in such a way that:
\begin{enumerate}
\item the curvature form of the smooth connection that we have defined on the boundary of each cube composing $C_{\varepsilon}$ is $\varepsilon$-close to the restriction of $F$ to that boundary;
\item such curvature forms patch together to give a well defined smooth $2$-form on $\partial C_{\varepsilon}$.
\end{enumerate} 
We recall that by "smooth" we always mean "which is the restriction of a smooth and compactly supported smooth form on $\r^3$". The proof of the previous facts is achieved through a sufficient condition to be the curvature of a weak connection on a Lipschitz principal $U(1)$-bundle over the boundary of a cube in $\r^3$ (Lemma \ref{L4.1}, see also \cite[Chapter IV, §20]{bott-tu}), standard convolution techniques and well known extension properties of Lipschitz maps (Lemma \ref{L4.2} and Lemma \ref{L4.3}).\\
In section 5, fixed any $\varepsilon>0$ and any good cube $Q$ from $C_{\varepsilon}$, we aim to replace $F$ by a smooth $2$-form on $Q$ which is, of course, $\varepsilon$-near to $F$ on $L^p(Q)$. The idea is to consider the smoothing $\phi$ of $F$ on $\partial Q$ (given by the previous section) and find a smooth $1$-form $\alpha$ on $\partial Q$ such that
\[\begin{cases}d\alpha=\phi\,\\
d^*\alpha=0,\end{cases} \mbox{ weakly on } \partial Q.\]
This is often called "gauge fixing" procedure. Then, we extend $\alpha$ to a smooth harmonic $1$-form $A$ on $Q$ which is also $C^{1,\alpha}$ up to $\partial Q$ for every $\alpha\in[0,1)$. Eventually, we replace $F$ with $dA$ on $Q$. The proofs of all the elliptic regularity results concerning the harmonic extension of a boundary datum in the interior of an open cube in $\r^n$ that we will use are collected in the Appendix A of the present paper, in order to keep the exposition cleaner.\\ 
Section 6 is dedicated to show that, fixed any $\varepsilon>0$, we can replace $F$ on every bad cube with the radial extension of the smoothing of $F$ on the boundary of that cube by remaining close to $F$ in $L^p\big(\{\mbox{union of the bad cubes in } C_{\varepsilon}\}\big)$. The $2$-form that we obtain by performing this replacement is smooth in the interior of every good cube up to the centre of the cube itself and Lipschitz on the union of all the bad cubes. Moreover, the singularity at the centre of every bad cube $Q$ composing $C_{\varepsilon}$ has degree $2\pi d$, where
\[d=\int_{\partial Q}i^*F.\]\\
In the first part of section 7, we just adjust the last technicalities to obtain a family of $2$-forms $\{F_{\varepsilon}\}_{\varepsilon\in E}$ which are defined on the whole unit ball $B$ and still satisfy the requirements of the statement. Namely, by radial extension starting from the boundary of the union of all the cubes in $C_{\varepsilon}$ we get a Lipschitz $2$-form on $B_2$ having the desired properties and by standard convolution and restriction to $B$ we eventually get the result.

\bigskip\noindent
The last part of the paper (subsections 7.1 and 7.2) is aimed to prove two remarkable consequences of Theorem \ref{T1.1}. The first one is the characterization of the $L^p$-vector fields with integer valued fluxes on the unit ball $B\subset\r^3$ when $p\ge 3/2$.
\begin{dfn}
We say that a vector field $X\in L^p(B;\r^3)$ has \text{integer valued fluxes} if 
\[\int_{\partial C_r(x_0)}\big(X\cdot\nu_{\partial C_r(x_0)}\big)\, d\H^2\in\z,\]
for every $x_0\in B$ and for $\L^1$-a.e. Lebesgue point of the function $X_{x_0}$ given by
\[X_{x_0}(r):=\int_{\partial C_r(x_0)}\big(X\cdot\nu_{\partial C_r(x_0)}\big)\, d\H^2, \qquad \mbox{ for every } r\in \bigg(0,\frac{1}{\sqrt{3}}\dist(x_0,\partial B)\bigg).\]
We denote the class of vector field with integer valued fluxes on $B$ by $L_{\z}^p(B)$.
\end{dfn}
\begin{dfn}
\label{definition vector fileds integer valued fluxes}
We say that a vector field $X\in L^p(B;\r^3)$ is \textbf{smooth away from integer point singularities} if $\star(2\pi X^{\flat})\in\Omega_{p,R}^2(B)$. We denote by $L_R^p(B)$ this class of vector fields on $B$. 
\end{dfn}
\begin{rem}
\label{R1.3}
We notice that there exists a linear isomorphisms between $\Omega_{p,\z}^2(B)$ and $L_{\z}^p(B)$ given by
\[F\mapsto (\star F)^{\sharp}.\]
Thus, these two spaces can always be identified and an analogous strong $L^p$-approximation result for vector fields in $L_{\z}^p(B)$ with elements of $L_R^p(B)$ follows trivially from Theorem \ref{T1.1}. 
\end{rem}
\noindent 
In Section 7.1, we prove the following Corollary \ref{C1.2} regarding this special class of vector fields on $B$. Up the knowledge of the author, a proof of this result never appeared before in literature.
\begin{cor}
\label{C1.2}
Let $p\ge 3/2$ and $X$ be any vector field in $L_{\z}^p(B)$. Then, $X$ is divergence free.
\end{cor}
\begin{rem}
By Corollary \ref{C1.2}, it holds that
\[L_{\z}^p(B)\subset\big\{\mbox{divergence free vector fields in } L^p(B;\r^3)\big\}, \qquad \mbox{ for every } p\ge 3/2.\]
Thus, since by Lemma \ref{LB.2} we get that
\[\big\{\mbox{divergence free vector fields in } L^p(B;\r^3)\big\}\subset L_{\z}^p(B), \qquad \mbox{ for every } p\ge 1,\]
we conclude that
\[L_{\z}^p(B)=\big\{\mbox{divergence free vector fields in } L^p(B;\r^3)\big\}, \qquad \mbox{ for every } p\ge 3/2\]
and we have thus obtained a characterization of the class $L_{\z}^p(B)$, whenever $p\ge 3/2$. 
\end{rem}
\noindent 
The second consequence of Theorem \ref{T1.1} that we will investigate and prove is the following decomposition theorem for vector fields in $L_{\z}^1(B)$.
\begin{thm}
\label{T1.2}
Let $X\in L_{\z}^1(B)$. Then, there exists an integer $1$-current $L\in\mathcal{R}_1(B)$ with finite mass and a $1$-cycle $C\in\mathcal{D}_1(B)$ such that:
\begin{enumerate}
\item $T_X=C+L$;
\item $L$ is a mass-minimizer on the class $\mathcal{R}_1(B)\cap\{T\in\mathcal{D}_1(B) \mbox{ s.t. } \partial T=\partial T_X\}$ and its mass is given by
\[\mathbb{M}(L)=\sup_{\substack{\varphi\in\mathcal{D}(B), \\ ||\nabla\varphi||_{L^{\infty}(B)}\le 1}}\left<\partial T_X,\varphi\right>=\sup_{\substack{\varphi\in\mathcal{D}(B), \\ ||\nabla\varphi||_{L^{\infty}(B)}\le 1}}\int_B(X\cdot\nabla\varphi)\, d\L^3.\]
\end{enumerate}
\end{thm}
\noindent
Theorem \ref{T1.2} can be regarded as a generalization of the result presented in \cite{brezis-coron-lieb} concerning existence of minimal connections for vector fields in $L_R^1(B)$. 

\section{Basics on Sobolev principal $U(1)$-bundles over Lipschitz submanifolds in $\r^n$}
\noindent
Let $M\subset\r^n$ be a compact Lipschitz $k$-submanifold in $\r^n$ such that $\H^k(M)<+\infty$. Since $M$ is a $k$-rectifiable subset of $\r^n$, the approximate tangent space $T_xM$ to $M$ exists on some $\H^k$-measurable subset $E_M\subset\r^n$ such that $\H^k(M\smallsetminus E_M)=0$. Whenever $T_xM$ exists, we denote by $\pi_x^M:\r^n\rightarrow T_xM$ the orthogonal projection from $\r^n$ to $T_xM$.
Given any $u\in C^{\infty}(M)$, we let the approximate gradient $\nabla u:M\rightarrow \r^n$ of $u$ on $M$ be defined as
\[\nabla u(x):=\begin{cases}\pi_x^M\big(\nabla\tilde u(x)\big) & \mbox{ on } E_M\\0 & \mbox{ on } M\smallsetminus E_M,\end{cases}\] 
where $\tilde u$ is any function in $C_c^{\infty}(\r^n)$ such that $u=\tilde u|_M$.\\ 
For every $u\in C^{\infty}(M)$, we let the $W^{1,p}$-Sobolev norm of $u$ be given by
\[||u||_{W^{1,p}(M)}=||u||_{L^p(M,\H^k\lfloor_M)}+||\nabla u||_{L^p(M,\H^k\lfloor_M)}.\]
\begin{dfn}[Sobolev spaces on Lipschitz compact submanifolds in $\r^n$]
Let $M$ be a Lipschitz compact $k$-submanifold in $\r^n$ such that $\H^k(M)<+\infty$ and $1\le p\le +\infty$.\\
We define the Sobolev space $W^{1,p}(M)$ as
\[W^{1,p}(M):=\overline{C^{\infty}(M)}^{||\cdot||_{W^{1,p}(M)}}.\]
For every $l\in\n$, we define the Sobolev space $W^{1,p}\big(M;\r^l\big)$ as the space of the functions $u:M\rightarrow\r^l$ such that each component of $u$ belong to $W^{1,p}(M)$.\\
If $N\subset\r^l$ is any isometrically embedded submanifold of $\r^l$, we define the Sobolev space $W^{1,p}(M;N)$ as 
\[W^{1,p}(M;N):=\Big\{u\in W^{1,p}\big(M,\r^l\big)\mbox{ s.t. } u(x)\in N \mbox{ for } \vol_g\mbox{-a.e. } x\in M\Big\}.\] 
\end{dfn}
\begin{dfn}[Lipschitz principal $U(1)$-bundles]
\label{D2.2}
Let $M$ be a compact Lipschitz $k$-submanifold in $\r^n$ such that $\H^k(M)<+\infty$.\\
A \textbf{Lipschitz principal} $U(1)$-\textbf{bundle} $P_{\mathcal{C},\rho}$ on $M$ is given by an relatively open covering $\mathcal{C}=\{U_{\alpha}\}_{\alpha\in A}$ of $M$ and a $U(1)$-valued function $\rho_{\alpha\beta}:U_{\alpha}\cap U_{\beta}\rightarrow U(1)$ on every non-empty intersection $U_{\alpha}\cap U_{\beta}$ such that:
\begin{enumerate}
\item $\rho_{\alpha\beta}\in W^{1,\infty}\big(U_{\alpha}\cap U_{\beta};U(1)\big)$ for every $\alpha,\beta\in A$ such that $U_{\alpha}\cap U_{\beta}\neq\emptyset$;
\item for every $\alpha,\beta,\gamma\in A$ such that $U_{\alpha}\cap U_{\beta}\cap U_{\gamma}\neq\emptyset$ it holds that $\rho_{\alpha\beta}\cdot \rho_{\beta\gamma}=\rho_{\alpha\gamma}$ as functions in $W^{1,p}\big(U_{\alpha}\cap U_{\beta}\cap U_{\gamma};U(1)\big)$. 
\end{enumerate}
We denote by $\P_{1,\infty}^{U(1)}(M)$ the set of all the Lipschitz principal $U(1)$-bundles on $M$.
\end{dfn}
\begin{rem}
Since the Lie group $U(1)$ is commutative, we also say that the elements of $\P_{1,\infty}^{U(1)}(M)$ are Lipschitz abelian bundles on $M$. 
\end{rem}
\noindent
One may wonder if the previous Definition \ref{D2.4} is well-posed. Indeed, although the inversion map $\cdot^{-1}:U(1)\rightarrow U(1)$ and the multiplication map $\cdot:U(1)\times U(1)\rightarrow U(1)$ are smooth on $U(1)$, a priori it is not completely obvious that $u^{-1}, u\cdot v\in W^{1,\infty}(M;U(1))$, for every couple of Sobolev maps $u,v\in W^{1,\infty}\big(M;U(1)\big)$. By However, the Gagliardo-Nirenberg interpolation inequality (see e.g. \cite[Comments to Chapter 9, C.]{brezis}) ensures that the compactness of $U(1)$ is enough to erase this kind of issues.
\begin{dfn}[Weak connections on Lipschitz abelian bundles]
\label{D2.3}
Let $M$ be a compact Lipschitz $k$-submanifold in $\r^n$ such that $\H^k(M)<+\infty$, $1\le p\le +\infty$ and $P_{\mathcal{C},\rho}\in\mathcal{P}_{1,\infty}^G(M)$.\\
A $\bm{W^{1,p}}$-\textbf{Sobolev connection} on $P_{\mathcal{C},\rho}$ is a family $iA=\{iA_j\}_{j\in J}$ of $\mathfrak{u}(1)$-valued local $1$-forms on $M$ such that:
\begin{enumerate}
\item $A_j\in W^{1,p}(U_j)$, for every $j\in J$;
\item the compatibility condition $iA_j=iA_i+\rho_{ij}^{-1}d\rho_{ij}$ is satisfied distributionally on $U_i\cap U_j$ for every $i,j\in J$ such that $U_i\cap U_j\neq\emptyset$. 
\end{enumerate}
Given any functional space $\mathscr{F}$ (e.g. $L^p$, $L^{p,\infty}$), by simply replacing $W^{1,p}$ with $\mathscr{F}$ in (1) we get the definition of $\mathscr{F}$-\textbf{connection} on $M$.\\
We say that an $L^1$-connection on $M$ is a \textbf{weak connection} on $M$. 
\end{dfn}
\begin{dfn}[Weak curvatures]
\label{D2.4}
Let $M$ be a compact Lipschitz $k$-submanifold in $\r^n$ such that $\H^k(M)<+\infty$, $1\le p\le +\infty$, $P_{\mathcal{C},\rho}\in\mathcal{P}_{1,\infty}^G(M)$ and $iA=\{iA_j\}_{j\in J}$ a weak connection on $P_{\mathcal{C},\rho}$.\\
If $dA_j\in L^p(U_j)$ for every $j\in J$, then we say that $iA$ admits a weak $L^p$-curvature and that the collection $F(iA):=\{idA_j\}_{j\in J}$ is the \textbf{weak} $L^p$-\textbf{curvature} of $A$.  
\end{dfn}
\begin{rem}
As a direct consequence of the compatibility condition (2) in Definition \ref{D2.3}, given a weak curvature $iA=\{iA_j\}_{j\in J}$ that admits a weak $L^p$-curvature then
\[dA_i=dA_j \qquad \mbox{ distributionally on } U_i\cap U_j,\]
on every non-empty intersection $U_i\cap U_j$.\\
Hence, the local $2$-forms $idA_j$ patch together to give a well-defined global $2$-form on $M$. We will denote such $2$-form again with $F(iA)$, with a slight abuse of notation.  
\end{rem}
\section{Choice of the cubic decompositions} 
\begin{figure}
	\centering
	\tdplotsetmaincoords{80}{130}
	\begin{tikzpicture}[scale=2, tdplot_main_coords]
	\coordinate (A) at (0.5,0.5,-0.5);
	\coordinate (B) at (0.5,0.5,0.5);
	\coordinate (C) at (0.5,-0.5,0.5);
	\coordinate (D) at (0.5,-0.5,-0.5);
	\coordinate (D) at (0.5,-0.5,-0.5);
	\coordinate (E) at (-0.5,-0.5,-0.5);
	\coordinate (F) at (-0.5,-0.5,0.5);
	\coordinate (G) at (0.5,-0.5,0.5);
	\coordinate (H) at (-0.5,0.5,-0.5);
	\coordinate (I) at (-0.5,0.5,0.5);
	\draw[color=black,->] (-1,0,0) -- (1,0,0) node[anchor=north east]{$x_1$};
	\draw[color=black,->] (0,-1,0) -- (0,1,0) node[anchor=north west]{$x_2$}; 
	\draw[color=black,->] (0,0,-1) -- (0,0,1) node[anchor=south]{$x_3$};
	\draw[color=black,thick,-] (A) -- (B);
	\draw[color=black,thick,-] (B) -- (C);
	\draw[color=black,thick,-] (C)-- (D);
	\draw[color=black,thick,-] (D) -- (A);
	\draw[fill=gray,opacity=0.4] (A) -- (B) -- (C) -- (D);
	\draw[color=black,thick,-] (D) -- (E);
	\draw[color=black,thick,-] (E) -- (F);
	\draw[color=black,thick,-] (F) -- (G);
	\draw[fill=gray,opacity=0.4] (D) -- (E) -- (F) -- (G);
	\draw[color=black,thick,-] (E) -- (H);
	\draw[color=black,thick,-] (H) -- (A);
	\draw[fill=gray,opacity=0.4] (D) -- (E) -- (H) -- (A);
	\draw[color=black,thick,-] (H) -- (I);
	\draw[color=black,thick,-] (I) -- (B);
	\draw[color=black,thick,-] (I) -- (F);
	\node at (G) [above left] {$\Gamma$};
	\end{tikzpicture}
	\caption{The picture represents the unit cube $C$ centred at the origin of $\r^3$. The part of $\partial C$ which is filled in gray is the Lipschitz surface $\Gamma$.} \label{fig1}
\end{figure}
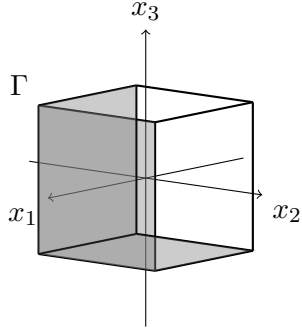
\noindent
We will need to consider the Lipschitz surface in $\r^3$ given by
\[\Gamma:=\partial C\smallsetminus\Big(\big\{x^3=1/2\big\}\cup\big\{x^2=1/2\big\}\cup\big\{x^1=-1/2\big\}\Big),\]
which is represented in Figure \ref{fig1} for the reader's convenience.\\
We set
\[\Gamma_{\varepsilon}(x_0):=\varepsilon\Gamma+x_0\]
\[C_\varepsilon^l(x_0):=C_\varepsilon(x_0)\cup\Gamma_\varepsilon(x_0),\]
and
\[C_\varepsilon^r(x_0):=C_\varepsilon(x_0)\cup\big(\partial C_{\varepsilon}(x_0)\smallsetminus\Gamma_\varepsilon(x_0)\big),\]
for every $x_0\in\r^3$ and for every $\varepsilon>0$.
Fix any $0<\varepsilon<1$ and consider the lattice $L_{\varepsilon}$ in $\r^3$ given by
\[L_{\varepsilon}:=\varepsilon\z^3+\bigg(\frac{\varepsilon}{2},\frac{\varepsilon}{2},\frac{\varepsilon}{2}\bigg).\]
We define $S_{\varepsilon}\subset\r^3$ as $S_{\varepsilon}:=L_{\varepsilon}\cap B_{1-3\varepsilon}(0)$ and we consider the following subsets of $B$:
\[C_{\varepsilon}^0:=\bigcup_{x_0\in S_{\varepsilon}}C_{\varepsilon}(x_0),\]
\[\partial C_{\varepsilon}^0:=\bigcup_{x_0\in S_{\varepsilon}}\partial C_{\varepsilon}(x_0).\]
We divide the points in $S_{\varepsilon}$ in the following 4 families:
\[\prescript{l}{}{S}_{\varepsilon}^1:=\big\{x_0\in S_{\varepsilon} \mbox{ s.t. } \partial C_{\varepsilon}(x_0)\smallsetminus\Gamma_{\varepsilon}(x_0) \mbox{ intersects the boundary of at least two cubes in } C_{\varepsilon}^0\big\},\]
\[\prescript{l}{}{S}_{\varepsilon}^2:=S_{\varepsilon}\smallsetminus \prescript{l}{}{S}_{\varepsilon}^1,\]
\[\prescript{r}{}{S}_{\varepsilon}^1:=\big\{x_0\in S_{\varepsilon} \mbox{ s.t. } \Gamma_{\varepsilon}(x_0) \mbox{ intersects to the boundary of at least two cubes in } C_{\varepsilon}^0\big\},\]
\[\prescript{r}{}{S}_{\varepsilon}^2:=S_{\varepsilon}\smallsetminus \prescript{r}{}{S}_{\varepsilon}^1.\]
Eventually, we define
\[C_{\varepsilon,a}:=C_{\varepsilon}^0+a\]
and
\[\partial C_{\varepsilon,a}:=\partial C_{\varepsilon}^0+a,\]
for every fixed $a\in\overline{B_{\varepsilon}(0)}$.
\begin{lem}
\label{L3.1}
Let $F\in\Omega_{p,\z}^2(B)$ and fix $\varepsilon\in (0,1)$. For every $a\in\overline{B_{\varepsilon}(0)}$, consider the piecewise-constant $2$-form $\prescript{l}{}{\bar F}_{\varepsilon,a}$ defined on $\overline{C_{\varepsilon,a}}$ by
\[\prescript{l}{}{\bar F}_{\varepsilon,a}(x):=\fint_{C_{\varepsilon}(x_0)+a}F(y)\, dy \qquad \mbox{ for every } x\in C_{\varepsilon}^l(x_0+a)\]
if $x_0\in\prescript{l}{}{S}_{\varepsilon}^1$ and 
\[\prescript{l}{}{\bar F}_{\varepsilon,a}(x):=\fint_{C_{\varepsilon}(x_0)+a}F(y)\, dy \qquad \mbox{ for every } x\in\overline{C_{\varepsilon}(x_0+a)}\]
otherwise. Analogously, define $\prescript{r}{}{\bar F}_{\varepsilon,a}$ on $\overline{C_{\varepsilon,a}}$ by
\[\prescript{r}{}{\bar F}_{\varepsilon,a}(x):=\fint_{C_{\varepsilon}(x_0)+a}F(y)\, dy \qquad \mbox{ for every } x\in C_{\varepsilon}^r(x_0+a)\]
if $x_0\in\prescript{r}{}{S}_{\varepsilon}^1$ and 
\[\prescript{r}{}{\bar F}_{\varepsilon,a}(x):=\fint_{C_{\varepsilon}(x_0)+a}F(y)\, dy. \qquad \mbox{ for every } x\in\overline{C_{\varepsilon}(x_0+a)}\]
otherwise. Then, there exist a subset $R\subset(0,1)$ such that $\L^1\big((0,1)\smallsetminus R\big)=0$ and a family of translations $\{a_{\varepsilon}\}_{\varepsilon\in R}\subset\overline{B_{\varepsilon}(0)}$ such that
\begin{equation}
\frac{1}{2\pi}\int_{\partial C_{\varepsilon}(x_0+a_{\varepsilon})} i_{\partial C_{\varepsilon}(x_0+a_{\varepsilon})}^*F\in\z, \qquad \mbox{ for every }\varepsilon\in R,
\end{equation}
and
\begin{align}\lim_{\varepsilon\rightarrow 0^+}\varepsilon||F-\prescript{l}{}{\bar F}_{\varepsilon,a_{\varepsilon}}||_{L^p(\partial C_{\varepsilon,a_{\varepsilon}})}^p=0,
\label{eq3.1}
\end{align}
\begin{align}\lim_{\varepsilon\rightarrow 0^+}\varepsilon||F-\prescript{r}{}{\bar F}_{\varepsilon,a_{\varepsilon}}||_{L^p(\partial C_{\varepsilon,a_{\varepsilon}})}^p=0.
\label{eq3.2}
\end{align}
Moreover, for such a choice, the following estimates hold
\begin{align}||F||_{L^p(\partial C_{\varepsilon,a_{\varepsilon}})}^p\le\frac{\kappa}{\varepsilon}||F||_{L^p(B)}^p+2^{p-1}||F- \prescript{l}{}{\bar F}_{\varepsilon,a_{\varepsilon}}||_{L^p(\partial C_{\varepsilon,a_{\varepsilon}})}^p,
\label{eq3.3}
\end{align}
\begin{align}
&||F||_{L^p(\partial C_{\varepsilon,a_{\varepsilon}})}^p\le\frac{\kappa}{\varepsilon}||F||_{L^p(B)}^p+2^{p-1}||F- \prescript{r}{}{\bar F}_{\varepsilon,a_{\varepsilon}}||_{L^p(\partial C_{\varepsilon,a_{\varepsilon}})}^p,
\label{eq3.4}
\end{align}
for some $\kappa>0$.
\begin{proof}
Denote by $E$ the following subset of $\r\times\r^3$
\[E:=\bigg\{(r,a)\in (0,1)\times\r^3 \mbox{ s.t. } a\in\overline{B_r(0)} \mbox{ and } \frac{1}{2\pi}\int_{\partial C_r(x_0+a)}i_{\partial C_r(x_0+a)}^*F\notin\z \mbox{ for some } x_0\in S_r\bigg\}.\]
For every fixed $r\in (0,1)$ and $a\in\r^3$ we define the following sections of $E$:
\[E_r:=\{a\in\r^3 \mbox{ s.t. } (r,a)\in E\},\]
\[E_a:=\{r\in(0,1) \mbox{ s.t. } (r,a)\in E\},\]
We claim that for $\L^1$-a.e. $r\in (0,1)$ it holds that $\L^3(E_r)=0$. This follows by Fubini's theorem, since
\[\int_0^1\L^3(E_r)\, dr=(\L^1\times\L^3)(E)=\int_{\r^3}\L^1(E_a)\, da\]
and by hypothesis we know that $\L^1(E_a)=0$ for every fixed $a\in\r^3$. Thus, it follows that there exists a subset $R\subset (0,1)$ such that $\L^1\big((0,1)\smallsetminus R\big)=0$ and for every fixed $\varepsilon\in R$ it holds that 
\[\frac{1}{2\pi}\int_{\partial C_{\varepsilon}(x_0+a)}i_{\partial C_{\varepsilon}(x_0+a)}^*F\in\z\]
for every $x_0\in S_{\varepsilon}$ and for $\L^3$-a.e. $a\in\overline{B_{\varepsilon}(0)}$.\\
In order to prove that we can provide a good family of translations $\{a_{\varepsilon}\}_{\varepsilon\in R}$ in $\overline{B_{\varepsilon}(0)}$ such that \ref{eq3.1} and \ref{eq3.3} hold, we first claim that  
\[\int_{\overline{B_{\varepsilon}(0)}}\int_{\partial C_{\varepsilon,a}}\left|F(x)-\prescript{l}{}{\bar F}_{\varepsilon,a}(x)\right|^p\ d\mathcal{H}^2(x)\, da=o\big(\varepsilon^2\big)\]
Indeed, let $\Sigma_{\varepsilon}^1(x_0):=\Gamma_{\varepsilon}(x_0)$ and $\Sigma_{\varepsilon}^2(x_0):=\partial C_{\varepsilon}(x_0)$, for every $x_0\in\r^3$ and $\varepsilon>0$. By Fubini's theorem, we get that
\begingroup
\allowdisplaybreaks
\begin{align*}
\int_{\overline{B_{\varepsilon}(0)}}\int_{\partial C_{\varepsilon,a}}\big|F-\prescript{l}{}{\bar F}_{\varepsilon,a}\big|^p\, d\mathcal{H}^2\, da&=\sum_{j=1}^2\bigintssss_{\overline{B_{\varepsilon}(0)}}\sum_{x_0\in S_{\varepsilon}^j}\bigintssss_{\Sigma_{\varepsilon}^j(x_0)+a}\left|F-\fint_{C_{\varepsilon}(x_0)+a}F\right|^p\, d\mathcal{H}^2\, da\\
&=\sum_{j=1}^2\bigintssss_{\overline{B_{\varepsilon}(0)}}\sum_{x_0\in S_{\varepsilon}^j}\bigintssss_{\Sigma_{\varepsilon}^j(0)}\left|F(\hspace{0.52mm}\cdot+x_0+a)-\fint_{C_{\varepsilon}(x_0)+a}F\right|^p\, d\mathcal{H}^2\, da\\
&=\sum_{j=1}^2\bigintssss_{\Sigma_{\varepsilon}^j(0)}\sum_{x_0\in S_{\varepsilon}^j}\bigintssss_{\overline{ B_{\varepsilon}(0)}}\left|F(\hspace{0.52mm}\cdot+x_0+a)-\fint_{C_{\varepsilon}(x_0)+a}F\right|^p\, da\, d\mathcal{H}^2\\
&=\sum_{j=1}^2\bigintssss_{\Sigma_{\varepsilon}^j(0)}\sum_{x_0\in S_{\varepsilon}^j}\bigintssss_{\overline{ B_{\varepsilon}(x_0)}}\left|F(x+y)-\fint_{C_{\varepsilon}(y)}F\right|^p\, dy\, d\mathcal{H}^2(x)\\
&=2^{p-1}\sum_{j=1}^2\bigintssss_{\Sigma_{\varepsilon}^j(0)}\Bigg(\sum_{x_0\in S_{\varepsilon}^j}\int_{\overline{B_{\varepsilon}(x_0)}}\big|F(x+y)-F(y)\big|^p\, dy\\
&\phantom{{}={}}+\sum_{x_0\in S_{\varepsilon}^j}\bigintssss_{\overline{ B_{\varepsilon}(x_0)}}\left|F(y)-\fint_{C_{\varepsilon}(y)}F\right|^p\, dy\Bigg)\, d\mathcal{H}^2(x)\\
&\le 2^{p-1}\sum_{j=1}^2\bigintssss_{\Sigma_{\varepsilon}^j(0)}\Bigg(\int_B\big|\tilde F(x+y)-\tilde F(y)\big|^p\, dy\\
&\phantom{{}={}}+\int_{\overline{B_{1-\varepsilon}(0)}}\left|F(y)-\fint_{C_{\varepsilon}(y)}F\right|^p\, dy\Bigg)\, d\mathcal{H}^2(x),
\end{align*}
where we have denoted by $\tilde F$ the extension of $F$ by $0$ outside $\overline{B_{1-\varepsilon}(0)}$.\\
\endgroup
Since $\tilde F\in\Omega_{L^p}^2(B)$, it holds that
\begin{align*}
\int_{B}\big|\tilde F(x+y)-\tilde F(y)\big|^p\, dy=\big|\big|\tilde F(x+\cdot)-\tilde F\big|\big|_{L^p(B)}^p=o(1) \qquad \mbox{ as } x\rightarrow 0.
\end{align*}
Moreover, fix any $y\in\overline{B_{1-\varepsilon}(0)}$ and, for every $\varepsilon\in R$, define the function  $\displaystyle{\rchi_{\varepsilon}:=\frac{1}{\varepsilon^3}\rchi_{C_1(y)}\bigg(\frac{\cdot}{\varepsilon}\bigg)}$ on $\r^3$. Notice that
\[\big(F\star\rchi_{\varepsilon}\big)(y):=\fint_{C_{\varepsilon}(y)}F, \qquad \mbox{ for every } y\in\overline{B_{1-\varepsilon}(0)}.\]
As $F\in\Omega_{L^p}^2(B)$ by hypothesis and since $\{\rchi_{\varepsilon}\}_{\varepsilon\in R}$ is an $L^1$-approximate identity, it is well known that
\[\big|\big|F-F\star\rchi_{\varepsilon}\big|\big|_{L^p\left(\overline{B_{1-\varepsilon}(0)}\right)}^p=\int_{\overline{B_{1-\varepsilon}(0)}}\left|F(y)-\fint_{C_{\varepsilon}(y)}F\right|^p\, dy=o(1), \qquad \mbox{ as } \varepsilon\rightarrow 0^+.\]
Hence, for every fixed $\delta>0$ there exists some $\varepsilon>0$ such that both
\[\int_{B}\big|\tilde F(x+y)-\tilde F(y)\big|^p\, dy<\frac{\delta}{2^p},\]
for every $x\in B$ such that $|x|<\varepsilon$, and 
\[\int_{\overline{B_{1-\varepsilon}(0)}}\left|F(y)-\fint_{C_{\varepsilon}(y)}F\right|^p\, dy<\frac{\delta}{2^p}.\]
This in turn implies that for every fixed $\delta>0$ there exists some $\varepsilon>0$ such that 
\[\frac{1}{\varepsilon^2}\int_{\overline{B_{\varepsilon}(0)}}\int_{\partial C_{\varepsilon,a}}\big|
F-\prescript{l}{}{\bar F}_{\varepsilon,a}\big|^p\, d\mathcal{H}^2\, da< \frac{1}{\varepsilon^2}\sum_{j=1}^2\int_{\Sigma_{\varepsilon}^j(0)}\delta\, d\H^2=9\delta\]
and our claim follows.\\
By the mean value theorem, we know that for every fixed $\varepsilon\in R$ there exists some non-negligible subset  $T_{\varepsilon}\subset\overline{B_{\varepsilon}(0)}$ such that both
\begin{equation}\big|\big|F-\prescript{l}{}{\bar F}_{\varepsilon,a}\big|\big|_{L^p(\partial C_{\varepsilon,a})}^p=\int_{\partial C_{\varepsilon,a}}\big|F-\prescript{l}{}{\bar F}_{\varepsilon,a}\big|^p\, d\mathcal{H}^2\le\frac{1}{\varepsilon^3}\int_{\overline{B_{\varepsilon}(0)}}\int_{\partial C_{\varepsilon,a}}\big|F-\prescript{l}{}{\bar F}_{\varepsilon,a}\big|^p\, d\mathcal{H}^2\, da.
\end{equation}
and 
\begin{equation}
\big|\big|\prescript{l}{}{\bar F}_{\varepsilon,a}\big|\big|_{L^p\left(\partial C_{\varepsilon,a}\right)}^p=\int_{\partial C_{\varepsilon,a}}\big|\prescript{l}{}{\bar F}_{\varepsilon,a}\big|^p\, d\mathcal{H}^2\le\frac{3}{\varepsilon}||F||_{L^p(B)}^p,
\label{eq3.5}
\end{equation}
hold for every $a\in T_{\varepsilon}$. By the previous claim, for every $a\in T_{\varepsilon}$ it holds that
\[\big|\big|F-\prescript{l}{}{\bar F}_{\varepsilon,a}\big|\big|_{L^p(\partial C_{\varepsilon,a})}^p=o\big(\varepsilon^{-1}\big)\]
which is equivalent to \ref{eq3.1}. Then, the estimate \ref{eq3.3} follows by triangle inequality and by \ref{eq3.5}.\\
Eventually, by applying exactly the same procedure to $\prescript{r}{}{\bar F}_{\varepsilon,a}$ we obtain that \ref{eq3.2} and \ref{eq3.4} holds for every $a\in T_{\varepsilon}'\subset T_{\varepsilon}$, where $T_{\varepsilon}'$ is a non-negligible subset of $\overline{B_{\varepsilon}(0)}$. Hence, since we have proved as a first fact that there exists $\tilde T_{\varepsilon}\subset \overline{B_{\varepsilon}(0)}$ such that $\L^3\big(\overline{B_{\varepsilon}(0)}\smallsetminus\tilde T_{\varepsilon}\big)=0$ and (3.1) holds for every $a\in\tilde T_{\varepsilon}$, we conclude by picking $a_{\varepsilon}\in T_{\varepsilon}'\cap\tilde T_{\varepsilon}$ (which is non-negligible and then non-empty).
\end{proof}
\end{lem}
\noindent
Fix any $F\in\Omega_{p,\z}^2(B)$. From now on, for every $\varepsilon\in R$ we simply denote by $C_{\varepsilon}$ the cubic decomposition $C_{\varepsilon,a_{\varepsilon}}$ provided by the previous Lemma \ref{L3.1}. Denoting by $N_{\varepsilon}:=\card(S_{\varepsilon})$, we  enumerate the $\varepsilon$-cubes composing $C_{\varepsilon}$ by the index $j\in\{1,...,N_{\varepsilon}\}$, so that
\[C_{\varepsilon}=\bigcup_{j=1}^{N_{\varepsilon}}C_{\varepsilon,j}.\]
For every fixed $j\in\{1,...,N_{\varepsilon}\}$, we say that the cube $C_{\varepsilon,j}$ is a good cube if 
\[\left|\frac{1}{2\pi}\int_{\partial C_{\varepsilon,j}}i_{\partial C_{\varepsilon,j}}^*F\right|<1\]
and a bad cube otherwise. We denote by $C_{\varepsilon}^g$ the interior part of following subset of $\r^3$:
\[\bigcup_{\substack{Q \subset C_{\varepsilon}\\ Q \text{ good cube}}}\overline{Q}.\]
Analogously, we define $C_{\varepsilon}^b$ by replacing "good" with "bad". Finally, we let $N_{\varepsilon}^g$ (resp. $N_{\varepsilon}^b$) be the number of good (resp. bad) cubes in $C_{\varepsilon}$.
\begin{lem}
\label{L3.2}
The volume of the union of the bad cubes vanishes as $\varepsilon\rightarrow 0^+$ in $R$, i.e.
\[\lim_{\varepsilon\rightarrow 0^+}\mathcal{L}^3\big(C_{\varepsilon}^b\big)=0\]
\begin{proof}
For every bad cube $C_{\varepsilon,j}$, it holds that
\[1\le\left|\frac{1}{2\pi}\int_{\partial C_{\varepsilon,j}}i_{\partial C_{\varepsilon,j}}^*F\right|\le\frac{1}{2\pi}\int_{\partial C_{\varepsilon,j}}\big|F\big|\, d\H^2.\]
Thus, by summing over all the bad cubes composing $C_{\varepsilon}$ 
\[(N_{\varepsilon}^b)^{1/p}\le\frac{1}{2\pi}\int_{\partial C_{\varepsilon}^b}\big|F\big|\, d\H^2\le\frac{1}{2\pi}||F||_{L^p\left(\partial C_{\varepsilon}^b\right)}\big(6\varepsilon^2\big)^{1/p'}.\]
Exploiting the estimate \ref{eq3.3} provided by Lemma \ref{L3.1}, we get
\begin{align*}
\L^3\big(C_{\varepsilon}^b\big)=N_{\varepsilon}^b\varepsilon^3\le\frac{6^{p-1}}{(2\pi)^p}\varepsilon^{2p+1}||F||_{L^p\left(\partial C_{\varepsilon}^b\right)}^p\le\Big(\kappa_1||F||_{L^p(B)}^p+ \kappa_2\varepsilon\big|\big|F-\prescript{l}{}{\bar F}_{\varepsilon}\big|\big|_{L^p(\partial C_{\varepsilon}^b)}^p\Big)\varepsilon^{2p}
\end{align*}
for some $\kappa_1,\kappa_2>0$ and the statement follows.
\end{proof}
\end{lem}
\noindent

\section{Smoothing on the boundary of the cubic decompositions}
\begin{lem}
\label{newlem}
Let $Q\subset\r^3$ be an open cube of side-length $l>0$ and centre $x_0=\big(x_0^1,x_0^2,x_0^3\big)$ in $\r^3$. Let 
\begin{align*}
U^+&:=\partial Q\cap\bigg\{x^3-x_0^3>-\frac{l}{4}\bigg\},\\
U^-&:=\partial Q\cap\bigg\{x^3-x_0^3<\frac{l}{4}\bigg\},\\
U^{\pm}&:=U^-\cap U^+.
\end{align*}
Then
\begin{enumerate}
	\item  if $p\in(1,+\infty)$, then for every $F\in\Omega_{L^p}^2(U^+)$ there exists $\alpha\in\Omega_{W^{1,p}}(U^+)$ such that $d\alpha=F$;
	\item the same statement as in (2) holds for $U^-$;
	\item if $p=1$, the same statements as in (1) and (2) hold up to requiring just $\alpha\in\Omega_{L^{2,\infty}}^1(U^+)$ and $\alpha\in\Omega_{L^{2,\infty}}^1(U^-)$ respectively;
\end{enumerate}
\begin{proof}
Since (2) and (3) are achieved by the same arguments used in (1), we just focus on (1) and (4).\\
Let $f:U^+\rightarrow B^2\subset\r^2$ be a bi-Lipschitz homeomorphism and consider the $2$-form on $B^2$ given by $\tilde F:=(f^{-1})^*F\in\Omega_{L^p}^2(B^2)$. Since $\tilde F$ is an $L^p$ top dimensional form on $B^2$, there exists $A\in\Omega_{W^{1,p}}^1(B^2)$ such that $dA=\tilde F$. The fact that $\alpha:=f^*A$ has the desired properties follows easily from the fact that, since $f$ is bi-Lipschitz, the operator $f^*$ commutes with weak differential. 
\end{proof}
\end{lem}
\begin{lem}
\label{L4.1}
Let $Q\subset\r^3$ be an open cube of side-length $l>0$ and centre $x_0=\big(x_0^1,x_0^2,x_0^3\big)$ in $\r^3$. Let $F\in\Omega_{L^p}^2\big(\partial Q\big)$ with $1\le p<+\infty$ be such that 
\[\frac{1}{2\pi}\int_{\partial Q}F\in\z.\]
Then:
\begin{enumerate}
\item if $p=1$, there exists $P_{\mathcal{C},\rho}\in\mathcal{P}_{1,\infty}^{U(1)}(\partial Q)$ and an $L^{2,\infty}$-connection $i\alpha$ on $P_{\mathcal{C},\rho}$ such that $F(i\alpha)=iF$.
\item if $1<p<+\infty$, there exist $P_{\mathcal{C},\rho}\in\mathcal{P}_{1,\infty}^{U(1)}(\partial Q)$ and a $W^{1,p}$-Sobolev connection $i\alpha$ on $P_{\mathcal{C},\rho}$ such that $F(i\alpha)=iF$.
\end{enumerate}
\begin{proof}
We start by proving (2). We consider the following open subsets of $\partial Q$:
\[U^+:=\partial Q\cap\bigg\{x^3-x_0^3>-\frac{l}{4}\bigg\},\]
\[U^-:=\partial Q\cap\bigg\{x^3-x_0^3<\frac{l}{4}\bigg\}.\]
Clearly $\mathcal{C}:=\big\{U^-,U^+\big\}$ is an open cover of $\partial Q$ made of contractible sets. We claim that, since $F$ is an $L^p$ top dimensional form on $\partial Q$, there exist $\alpha^-\in\Omega_{W^{1,p}}^1\big(U^-\big)$ and $\alpha^+\in\Omega_{W^{1,p}}^1\big(U^+\big)$ such that
\begin{align}
\label{eq4.1}
d\alpha^-=i_{U^-}^*F \qquad \mbox{ on } U^-
\end{align}
and
\begin{align}
\label{eq4.2}
d\alpha^+=i_{U^+}^*F \qquad \mbox{ on } U^+.
\end{align}
We define the following closed subsets of $\partial Q$:
\[\partial Q^+:=\partial Q\cap\big\{x^3-x_0^3\ge 0\big\},\]
\[\partial Q^-:=\partial Q\cap\big\{x^3-x_0^3\le 0\big\},\]
\[\partial Q^{\pm}:=\partial Q^+\cap\partial Q^-.\]
We define the projection $\pi:U^{\pm}\rightarrow\partial Q^{\pm}$ given by $\pi(x^1,x^2,x^3):=(x^1,x^2,x_0^3)$.
Thus, the $1$-form $\omega:=\pi^*\gamma$ is a generator of $H_{dR}^1\big(U^{\pm}\big)$, where $\gamma\in\Omega_{W^{1,\infty}}^1(\partial Q^{\pm})$ is such that
\begin{align*}
\int_{\partial Q^{\pm}}\gamma=2\pi.
\end{align*}
Since $\alpha^+-\alpha^-$ is a closed $1$-form on $U^{\pm}$ and $H_{dR}^1\big(U^{\pm}\big)\cong\z$, without losing generality (i.e. possibly absorbing exact terms either in $\alpha^+$ or in $\alpha-$ and keeping in mind the procedure that we have already used in Lemma \ref{newlem} in order to perform the computations) we can assume that
\[\alpha^+-\alpha^-=a\omega \qquad \mbox{ for some } a\in\r.\]
Then, by Stokes theorem, it holds that
\begin{align*}
\z\ni n:=\frac{1}{2\pi}\int_{\partial Q}i_{\partial Q}^*F&=\frac{1}{2\pi}\int_{\partial Q-}i_{U^-}^*F+\frac{1}{2\pi}\int_{\partial Q^+}i_{U^+}^*F=\frac{1}{2\pi}\int_{\partial Q^-}d\alpha^-+\frac{1}{2\pi}\int_{\partial Q+}d\alpha^+=\\
&=\frac{1}{2\pi}\int_{\partial Q^{\pm}}i_{\partial Q^{\pm}}^*\alpha^+-\frac{1}{2\pi}\int_{\partial Q^{\pm}}i_{\partial Q^{\pm}}^*\alpha^-=\frac{a}{2\pi}\int_{\partial Q^{\pm}}i_{\partial Q^{\pm}}^*\omega=\frac{a}{2\pi}\int_{\partial Q^{\pm}}\gamma=a.
\end{align*}
We analyze separately the following two cases.
\begin{enumerate}
\item If $n=0$, then $\alpha^+$ and $\alpha^-$ patch together to give a globally well defined $1$-form $\alpha$ on $\partial Q$ such that $d\alpha=F$. This implies that $iF$ is the $L^p$-curvature of the $W^{1,p}$-connection $i\alpha$ on the trivial $U(1)$-principal bundle over $\partial Q$.
\item If $n\neq 0$, then we consider the smooth function $g:S^1\rightarrow U(1)$ given by
\[g(\cos\theta,\sin\theta):=e^{in\theta}, \qquad \mbox{ for every } \theta\in[0,2\pi]\]
and we define the Lipschitz transition function $\rho^{\pm}:U^{\pm}\rightarrow U(1)$ given by $\rho^{\pm}:=h^*g$. Then, $P_{\mathcal{C},\rho^{\pm}}$ is a Lipschitz $U(1)$-principal bundle on $\partial Q$ and a direct computation shows that $i\alpha=\{i\alpha^-,i\alpha^+\}$ is a $W^{1,p}$-connection on $P_{\mathcal{C},\rho^{\pm}}$ such that $F(i\alpha)=iF$. 
\end{enumerate}
This completes the proof of (2). For was concerns (1), the proof is identical except for the fact that the best regularity that we can guarantee for the local 1-forms $\alpha^+$ and $\alpha^-$ solving distributionally the differential equations \eqref{eq4.1} and \eqref{eq4.2} is given by $\alpha^+\in\Omega_{L^{2,\infty}}^1(U^+)$ and $\alpha^-\in\Omega_{L^{2,\infty}}^1(U^-)$.
\end{proof}
\end{lem}
\begin{lem}
\label{L4.2}
Let $Q\subset\r^3$ be an open cube of side-length $l>0$ and centre $x_0=\big(x_0^1,x_0^2,x_0^3\big)$ in $\r^3$. Consider a Lipschitz principal $U(1)$-bundle $P_{\mathcal{C},\rho}$ on $\partial Q$ given by the finite open covering $\mathcal{C}=\{U_j\}_{j\in J}$ and the family $\rho$ of transition functions. Eventually, pick any $L^p$-connection $i\alpha$ over $P_{\mathcal{C},\rho}$ that admits a weak $L^p$-curvature $F(i\alpha)$ for some $1\le p<+\infty$.\\
Then, for every $\varepsilon>0$ there exists a Lipschitz connection $i\alpha_{\varepsilon}$ over $P_{\mathcal{C},\rho}$ such that
\[\left|\left|F(i\alpha_{\varepsilon})-F(i\alpha)\right|\right|_{L^p(\partial Q)}<\varepsilon.\]
\begin{proof}
Consider a Lipschitz partition of unity $\{\lambda_j\}_{j\in J}$ on $\partial Q$ subordinated to the finite open covering $\mathcal{C}$ and set	
\[M:=\max_{j\in J}||\nabla\lambda_j||_{L^{\infty}(\partial Q)}.\]
Fix any $\varepsilon>0$. By standard convolution, for every $j\in J$ we get a Lipschitz $1$-form $\tilde\alpha_{\varepsilon}^j$ on $U_j$ such that 
\[\left|\left|\tilde\alpha_{\varepsilon}^j-\alpha^j\right|\right|_{L^p(U_j)}<\frac{\varepsilon}{2\card(J)^2M}\]	
and
\[\left|\left|d\tilde\alpha_{\varepsilon}^j-d\alpha^j\right|\right|_{L^p(U_j)}<\frac{\varepsilon}{2\card(J)^2}.\]		
However, the family $\big\{i\tilde\alpha_{\varepsilon}^j\big\}_{j\in J}$ doesn't necessarily define a connection on $P_{\mathcal{C},\rho}$. In order to solve this issue, for every fixed $j\in J$ we define 
\[i\alpha_{\varepsilon}^j:=\sum_{k\in J}\beta_{\varepsilon}^{jk},\]
where $\beta_{\varepsilon}^{jk}$ is the Lipschitz $1$-form on $U_j$ given by
\[\beta_{\varepsilon}^{jk}:=\begin{cases}\lambda_k\big(i\tilde\alpha_{\varepsilon}^k+\rho_{kj}^{-1}d\rho_{kj}\big), & \mbox{ on } U_j\cap U_k\\ 0 & \mbox{ on } U_j\smallsetminus U_k,\end{cases} \qquad \mbox{ for every } k\in J.\]
First, we check that the family $i\alpha_{\varepsilon}=\big\{i\alpha_{\varepsilon}^j\big\}_{j\in J}$ satisfies the compatibility condition
\begin{align}
\label{eq4.3}
i\alpha_{\varepsilon}^j=i\alpha_{\varepsilon}^k+\rho_{kj}^{-1}d\rho_{kj} \qquad \mbox{ on } U_j\cap U_k,
\end{align}
for every $j,k\in J$ such that $U_j\cap U_k\neq\emptyset$. To this purpose, we fix any $j,k\in J$ such that $U_j\cap U_k\neq\emptyset$ and we notice that on $U_j\cap U_k$ it holds that
\begin{align*}i\alpha_{\varepsilon}^k+\rho_{kj}^{-1}d\rho_{kj}&=i\alpha_{\varepsilon}^k+\rho_{jk}d\rho_{kj}=\sum_{l\in J}\beta_{\varepsilon}^{kl}+\rho_{jk}d\rho_{kj}
\end{align*}
By differentiating the cocycle condition
\[\rho_{lj}=\rho_{lk}\rho_{kj}\]
we get that
\[d\rho_{lj}=d\rho_{ik}\rho_{kj}+\rho_{lk}d\rho_{kj}.\]
Multiplying both sides by $\rho_{jl}$ and exploiting the commutativity of the product in $U(1)$, we get
\[\rho_{jl}d\rho_{lj}=\rho_{jl}\rho_{kj}d\rho_{lk}+\rho_{jl}\rho_{lk}d\rho_{kj}=\rho_{kl}d\rho_{lk}+\rho_{jk}d\rho_{kj}.\]
By definition of $\beta_{\varepsilon}^{kl}$ and of $\alpha_{\varepsilon}^j$, we get
\[\sum_{l\in J}\beta_{\varepsilon}^{kl}=\sum_{j\in J}\beta_{\varepsilon}^{jl}-\rho_{jk}d\rho_{kj}=i\alpha_{\varepsilon}^j-\rho_{jk}d\rho_{kj}.\]
Thus, the compatibility condition \eqref{eq4.3} is proved.\\
Finally, we claim that
\[\left|\left|d\alpha_{\varepsilon}^j-d\alpha^j\right|\right|_{L^p(U_j)}<\frac{\varepsilon}{\card(J)}, \qquad \mbox{ for every } j\in J.\]
Indeed, fix any $j\in J$ and notice that
\begin{align*}d\big(\lambda_k(i\tilde\alpha_{\varepsilon}^k+\rho_{kj}^{-1}d\rho_{kj})\big)&=d\lambda_k\wedge(i\tilde\alpha_{\varepsilon}^k+\rho_{kj}^{-1}d\rho_{kj})+i\lambda_kd\tilde\alpha_{\varepsilon}^k=id\lambda_k\wedge(\tilde\alpha_{\varepsilon}^k+\alpha^j-\alpha^k)+i\lambda_kd\tilde\alpha_{\varepsilon}^k\\
&=i\big(d\lambda_k\wedge(\tilde\alpha_{\varepsilon}^k-\alpha^k)+\lambda_kd\tilde\alpha_{\varepsilon}^k+d\lambda_k\wedge\alpha^j\big).
\end{align*}
Since
\[\sum_{k\in J}d\lambda_k=d\bigg(\sum_{k\in J}\lambda_k\bigg)=0,\]
it follows that
\[d\alpha_{\varepsilon}^j=\sum_{k\in J}\gamma_{\varepsilon}^{jk},\]
where
\[\gamma_{\varepsilon}^{jk}:=\begin{cases}d\lambda_k\wedge(\tilde\alpha_{\varepsilon}^k-\alpha^k)+\lambda_kd\tilde\alpha_{\varepsilon}^k, & \mbox{ on } U_j\cap U_k\\ 0 & \mbox{ on } U_j\smallsetminus U_k,\end{cases} \qquad \mbox{ for every } k\in J.\]
Hence,
\begin{align*}||d\alpha_{\varepsilon}^j-d\alpha^j||_{L^p(U_j)}&=\bigg|\bigg|d\alpha_{\varepsilon}^j-\sum_{k\in J}\lambda_k d\alpha^j\bigg|\bigg|_{L^p(U_j)}\le\sum_{k\in J}||\gamma_{\varepsilon}^{jk}-\lambda_kd\alpha^j||_{L^p(U_j)}\\
&=\sum_{k\in J}||\gamma_{\varepsilon}^{jk}-\lambda_kd\alpha^j||_{L^p(U_j\cap U_k)}=\sum_{k\in J}||d\lambda_k\wedge(\tilde\alpha_{\varepsilon}^k-\alpha^k)+\lambda_k(d\tilde\alpha_{\varepsilon}^k-d\alpha^j)||_{L^p(U_j\cap U_k)}\\	
&\le\sum_{k\in J}||d\lambda_k\wedge(\tilde\alpha_{\varepsilon}^k-\alpha^k)||_{L^p(U_j\cap U_k)}+\sum_{k\in J}||\lambda_k(d\tilde\alpha_{\varepsilon}^k-d\alpha^k)||_{L^p(U_j\cap U_k)}\\
&\le\sum_{k\in J}M||\tilde\alpha_{\varepsilon}^k-\alpha^k)||_{L^p(U_k)}+\sum_{k\in J}||d\tilde\alpha_{\varepsilon}^k-d\alpha^k||_{L^p(U_k)}<\frac{\varepsilon}{\card(J)}.
\end{align*}
Since $\{U_j\}_{j\in J}$ is a cover of $\partial Q$ and $F(i\alpha_{\varepsilon})=d\alpha_{\varepsilon}^j$, $F(i\alpha)=d\alpha^j$ on each $U_j$, we obtain that
\[\left|\left|F(i\alpha_{\varepsilon})-F(i\alpha)\right|\right|_{L^p(\partial Q)}\le\sum_{j\in J}\left|\left|F(i\alpha_{\varepsilon})-F(i\alpha)\right|\right|_{L^p(U_j)}=\sum_{j\in J}||d\alpha_{\varepsilon}^j-d\alpha^j||_{L^p(U_j)}\]
and the statement follows. 
\end{proof}
\end{lem}
\begin{lem}
\label{L4.3}
Let $Q\subset\r^3$ be any open cube of side-length $l>0$ and centre $x_0=\big(x_0^1,x_0^2,x_0^3\big)$ in $\r^3$ and $1\le p<+\infty$. Then, it holds that
\[\Lip\big(U^+;\r^3\big)\subset\overline{C^{\infty}\big(U^+;\r^3\big)}^{W^{1,p}}\]
and
\[\Lip\big(U^-;\r^3\big)\subset\overline{C^{\infty}\big(U^-;\r^3\big)}^{W^{1,p}}\]
\begin{proof}
Since the proof is identical for $U^-$, we just focus on the case of $U^+$. Let $f$ be any Lipschitz vector field on $U^+$. First of all, we extend $f$ to the open neighbourhood of $U^+$ given by
\[\Omega:=\{x\in\r^3 \mbox{ s.t. } \dist(x,U^+)<l/4\}\cap\{x^3-x_0^3>-l/4\}.\]
We remark that, in the following part of the proof, the "radial extension of $f$ at the point $x$ with centre $z$" will be given by $f(\pi_{z}(x))$, where $\pi_{z}:\r^3\rightarrow\partial C_l(z)$ is the projection map
\[\pi_{z}(x):=\frac{l}{2}\frac{x-z}{||x-z||_{*}}+z, \qquad \mbox{ for every } x\in\r^3\smallsetminus\{z\}.\]
We pick a radial extension of the datum $f$ on $U^+$, switching the centre basing on the position of the point to which we aim to achieve the extension. Namely, we surround $Q$ by copies of itself so that any point $x\in\Omega\smallsetminus\partial Q$ lies in either $Q$ or in one of these cubes. If $x\in Q$, then we extend radially with centre $x_0$ the datum $f\in\partial Q$ to $x$. If $x$ belongs to some cube surrounding $Q$, say $Q_x$ with centre $c_x$, then we extend radially with centre $c_x$ the datum $f\in\partial Q$ to $x$. In this way, we obtain a well-defined Lipschitz function $\hat f:\Omega\rightarrow\r^3$.\\
By McShane's lemma (see e.g. \cite[Lemma 7.3]{maggi}), there exists a Lipschitz extension $\tilde f:\r^3\rightarrow\r^3$ of $\hat f$ such that $\Lip\big(\tilde f\big)\le\sqrt{3}\Lip(\hat f,\Omega)$. Let $\{\rho_{\delta}\}_{0<\delta<1}\subset C_c^{\infty}(\r^3)$ be a regularizing kernel and define 
\[\tilde f_{\delta}:=\tilde f\star\rho_{\delta}, \qquad \mbox{ for every } 0<\delta<1.\]
Since $\tilde f$ is continuous and $U^+$ is a compact subset of $\r^3$, then we know that
\begin{align}
\label{eq4.4}
f_{\delta}:=\tilde f_{\delta}|_{\partial Q}\xrightarrow{\delta\rightarrow 0^+}\tilde f|_{\partial Q}=f \qquad \mbox{ uniformly.}
\end{align}
Moreover, we have that 
\begin{align*}
||\nabla f_{\delta}||_{L^{\infty}(\partial Q;\r^3)}&\le\Lip(f_{\delta},\partial Q)\le \Lip(\tilde f,\partial Q+\delta B)\\
&\le\Lip(\tilde f)\le\sqrt{3}\Lip(\hat f,\Omega)<+\infty, \qquad \mbox{ for every } 0<\delta<1,
\end{align*}
and $\nabla f_{\delta}(x)\xrightarrow{\delta\rightarrow 0^+}\nabla f(x)$ for $\H^2$-a.e. $x\in \partial Q$. Hence, by dominated convergence theorem, we conclude that $\nabla f_{\delta}\xrightarrow{\delta\rightarrow 0^+}\nabla f(x)$ strongly in $L^p\big(\partial Q,\H^2;\r^3\big)$. This and \ref{eq4.4} allow to conclude that $f_{\delta}\xrightarrow{\delta\rightarrow 0^+} f$ strongly in $W^{1,p}\big(\partial Q,\H^2;\r^3\big)$ and the statement follows.
\end{proof}
\end{lem}

\noindent 
We recall that by a "smooth connection" over a Lipschitz principal $U(1)$-bundle on the boundary of a cube we mean a Lipschitz connection $i\alpha$ on that bundle such that every local representation of $i\alpha$ is a restriction to its domain of a smooth and compactly supported $1$-form on $\r^3$.\\
The following corollary follows easily by Lemma \ref{L4.3} and by Lemma \ref{L4.2}, applying the same procedure that we have used in the proof of Lemma \ref{L4.2}.
\begin{cor}
\label{C4.1} 
Let $Q\subset\r^3$ be an open cube of side-length $l>0$ and centre $x_0=\big(x_0^1,x_0^2,x_0^3\big)$ in $\r^3$. Consider the Lipschitz principal $U(1)$-bundle $P_{\mathcal{C},\rho^{\pm}}$ on $\partial Q$ given by the open covering $\mathcal{C}=\{U^+,U^-\}$ and the transition function $\rho^{\pm}$ as it is defined in the proof of Lemma \ref{L4.1}. Eventually, consider any $L^p$-connection $i\alpha$ over $P_{\mathcal{C},\rho}$ that admits a weak $L^p$-curvature $F(i\alpha)$ for some $1\le p<+\infty$.\\
Then, for every $\varepsilon>0$ there exists a smooth connection $i\alpha_{\varepsilon}$ over $P_{\mathcal{C},\rho}$ such that
\[\left|\left|F(i\alpha_{\varepsilon})-F(i\alpha)\right|\right|_{L^p(\partial Q)}<\varepsilon.\]
\end{cor}
\begin{rem}
\label{R4.1}
Fix any $\varepsilon\in R$ and consider the cubic decomposition $C_{\varepsilon}$. While performing the proofs of the previous lemmata, we have been careful to always pick the regularizing kernels and the partitions of unity coherently. Therefore, the curvatures of the smooth $\varepsilon$-approximations of the $L^p$-connections on the boundary of the cubes composing $C_{\varepsilon}$ patch together to give a well defined smooth $2$-form on $\partial C_{\varepsilon}$.
\end{rem}

\noindent
Fix any $\varepsilon\in R$. We conclude this section by applying the previous lemmata and corollaries in order to perform the required smoothing on the boundary of the cubic decomposition $C_{\varepsilon}$.\\
By Lemma \ref{L4.1}, for every open cube $C_{\varepsilon,j}$ composing $C_{\varepsilon}$ we fix a Lipschitz principal $U(1)$-bundle $P_{\mathcal{C}_j\rho_j^{\pm}}$ on $\partial C_{\varepsilon,j}$ and a $L^p$-connection $i\alpha^j$ on $P_{\mathcal{C}_j,\rho_j^{\pm}}$ such that
\[F(i\alpha^j)=i_{\partial C_{\varepsilon,j}}^*F\]
By Corollary \ref{C4.1} and Remark \ref{R4.1}, for every $j\in N_{\varepsilon}$ we can select a smooth connection $i\hat\alpha_{\varepsilon}^j$ on $P_{\mathcal{C}_j,\rho_j^{\pm}}$ such that
\begin{align}
\label{eq4.5}
\left|\left|F\big(i\hat\alpha_{\varepsilon}^j\big)-F\big(i\alpha^j\big)\right|\right|_{L^p(\partial C_{\varepsilon,j})}<\varepsilon.
\end{align}
and the $2$-forms $-iF\big(i\hat\alpha_{\varepsilon}^j\big)$ patch together to give a well defined smooth function $\phi_{\varepsilon}$ on $\partial C_{\varepsilon}$.\\
For every $j\in N_{\varepsilon}$, we denote by $\phi_{\varepsilon}^j$ the restriction to $\partial C_{\varepsilon,j}$ of $\phi_{\varepsilon}$, namely, the $2$-form $-iF\big(i\hat\alpha_{\varepsilon}^j\big)$.
Moreover, for every good cube $C_{\varepsilon,j}$, by the proof of Lemma \ref{L4.1} we know that the principal $U(1)$-bundle that we have defined on $\partial C_{\varepsilon,j}$ is trivial. Thus, the local representations for the connection $i\hat\alpha_{\varepsilon}^j$ patch together to given a globally well defined $1$-form on $\partial C_{\varepsilon,j}$. With a slight abuse of notation, we indicate this $1$-form again with $i\hat\alpha_{\varepsilon}^j$. Moreover, we denote by $\hat A_{\varepsilon}^j\in\Omega_c^1(\r^3)$ the smooth and compactly supported $1$-form such that $\hat\alpha_{\varepsilon}^j=i_{\partial C_{\varepsilon,j}}^*\hat A_{\varepsilon}^j$. 
\noindent

\section{Harmonic approximation on the good cubes}
\noindent 
Let $C_{\varepsilon,j}$ be any good cube in the cubic decomposition $C_{\varepsilon}$ of size $\varepsilon\in R$. Let $\bar F_{\varepsilon}^j$ be the $2$-form on $\r^3$ defined by
\[\bar F_{\varepsilon}^j:=\fint_{C_{\varepsilon,j}}F \qquad \mbox{ on } \r^3.\]
Since $\bar F_{\varepsilon}^j$ is a constant $2$-form on $\r^3$, it holds that $d\bar F_{\varepsilon}^j=0$. As $\r^3$ is contractible, there exists $A\in\Omega^1\big(\r^3\big)$ such that 
\[dA=\bar F_{\varepsilon}^j \qquad \mbox{ on } \r^3.\]
Let $\bar A\in\Omega^1\big(\r^3\big)$ be a smooth $1$-form on $\r^3$ such that $\bar A\equiv A$ on $\overline{C_{\varepsilon,j}+B}$ and $\bar A\equiv 0$ on $\r^3\smallsetminus\big(\overline{C_{\varepsilon,j}+2B}\big)$ so that $d^*\bar A\in C_c^{\infty}\big(\r^3\big)$, Thus, by standard elliptic theory, there exists a smooth solution $\varphi\in C^{\infty}\big(\r^3\big)$ for the Poisson equation
\[d^*du=\Delta u=-d^*\bar A \qquad \mbox{ on } \r^3.\]
We define $\bar A_{\varepsilon}^j:=\bar A+d\varphi\in\Omega^1\big(\r^3\big)$ and we notice that $\bar A_{\varepsilon}^j$ satisfies the following system of differential equations:
\[\begin{cases}d\bar A_{\varepsilon}^j=\bar F_{\varepsilon}^j,\\d^*\bar A_{\varepsilon}^j=0,\end{cases} \qquad \mbox{ on } C_{\varepsilon,j}.\]
We let $\bar\alpha_{\varepsilon}^j:=i_{\partial C_{\varepsilon,j}}^*\bar A_{\varepsilon}^j\in\Omega_{W^{1,\infty}}^1\big(\partial C_{\varepsilon,j}\big)$. Clearly, $\bar\alpha_{\varepsilon}^j$ is such that:
\[\begin{cases}d\bar\alpha_{\varepsilon}^j=i_{\partial C_{\varepsilon,j}}^*\bar F_{\varepsilon}^j,\\d^*\bar\alpha_{\varepsilon}^j=0,\end{cases} \qquad \mbox{ weakly on } \partial C_{\varepsilon,j}.\]
Next, we want to find a suitable $1$-form $\tilde\alpha_{\varepsilon}^j\in\Omega_{W^{1,\infty}}^1\big(\partial C_{\varepsilon,j}\big)$ having the following properties:
\begin{align}\begin{cases}d\tilde\alpha_{\varepsilon}^j=\phi_{\varepsilon}^j,\\d^*\tilde\alpha_{\varepsilon}^j=0,\end{cases} \qquad \mbox{ weakly on } \partial C_{\varepsilon,j}.\end{align}
Since $d^*\hat A_{\varepsilon}^j\in C_c^{\infty}\big(\r^3\big)$, by standard elliptic theory there exists a smooth solution $f$ of the following PDE:
\[\Delta u=-d^*\hat A_{\varepsilon}^j \qquad \mbox{ on } \r^3.\]
We claim that the $1$-form $\tilde\alpha_{\varepsilon}^j:=\hat\alpha_{\varepsilon}^j+i_{\partial C_{\varepsilon,j}}^*df\in\Omega_{W^{1,\infty}}^1\big(\partial C_{\varepsilon,j}\big)$ satisfies (5.1). Indeed
\[d\tilde\alpha_{\varepsilon}^j=d\big(\hat\alpha_{\varepsilon}^j+i_{\partial C_{\varepsilon,j}}^*df\big)=\phi_{\varepsilon}^j+i_{\partial C_{\varepsilon,j}}^*(d^2f)=\phi_{\varepsilon}^j\]
and
\[d^*\tilde\alpha_{\varepsilon}^j=d^*\big(\hat\alpha_{\varepsilon}^j+i_{\partial C_{\varepsilon,j}}^*df\big)=i_{\partial C_{\varepsilon,j}}^*\big(d^*df+d^*\hat A_{\varepsilon}^j\big)=i_{\partial C_{\varepsilon,j}}^*\big(\Delta f+d^*\hat A_{\varepsilon}^j\big)=0.\]
Now, we consider the harmonic extension (see Definition \ref{DA.1}) of the boundary datum $\tilde\alpha_{\varepsilon}^j\in\Omega_{W^{1,\infty}}^1\big(\partial C_{\varepsilon,j}\big)$ to the open cube $C_{\varepsilon,j}$, namely, the solution $\tilde A_{\varepsilon}^j$ of the differential problem
\begin{align*}\begin{cases}\Delta A=0&\mbox{ on } C_{\varepsilon,j}\\i_{\partial C_{\varepsilon,j}}^*A=\tilde\alpha_{\varepsilon}^j &\mbox{ on } \partial C_{\varepsilon,j}.\end{cases}\end{align*}
By Lemma \ref{LA.3}, $\tilde A_{\varepsilon}^j\in\Omega_{W^{2,p}}^1\big(C_{\varepsilon,j}\big)\cap\Omega_{C^{1,\alpha}}^1\big(\overline{C_{\varepsilon,j}}\big)$, for every $\alpha\in[0,1)$. Moreover $i_{\partial C_{\varepsilon,j}}^*\big(d^*\tilde A_{\varepsilon}^j\big)=0$ on $\partial C_{\varepsilon,j}$.\\
We define the $1$-form $A_{\varepsilon}^j\in\Omega_{W^{2,p}}^1\big(C_{\varepsilon,j}\big)\cap\Omega_{C^{1,\alpha}}^1\big(\overline{C_{\varepsilon,j}}\big)$ by 
\[A_{\varepsilon}^j:=\tilde A_{\varepsilon}^j-\bar A_{\varepsilon}^j\]
and the Lipschitz $1$-form $\alpha_{\varepsilon}^j\in\Omega_{W^{1,\infty}}^1\big(C_{\varepsilon,j}\big)$ by
\[\alpha_{\varepsilon}^j:=i_{\partial C_{\varepsilon,j}}^*A_{\varepsilon}^j.\]
Notice that $\alpha_{\varepsilon}^j$ satisfies the conditions
\begin{align*}\begin{cases}d\alpha_{\varepsilon}^j=\phi_{\varepsilon}^j-i_{\partial C_{\varepsilon,j}}^*\bar F_{\varepsilon}^j\\d^*\alpha_{\varepsilon}^j=0\end{cases}\mbox{ on } \partial C_{\varepsilon,j}.\end{align*}
and $A_{\varepsilon}^j$ is a solution of the differential problem 
\begin{align}\begin{cases}\Delta A=0&\mbox{ on } C_{\varepsilon,j}\\i_{\partial C_{\varepsilon,j}}^*A=\tilde\alpha_{\varepsilon}^j-\bar\alpha_{\varepsilon}^j &\mbox{ on } \partial C_{\varepsilon,j}.\end{cases}\end{align}
Again, by Lemma \ref{LA.3}, $A_{\varepsilon}^j$ is the unique solution of the previous system, $i_{\partial C_{\varepsilon,j}}^*\big(d^*A_{\varepsilon}^j\big)=0$ and the estimate 
\begin{align}\big|\big|dA_{\varepsilon}^j\big|\big|_{L^p(C_{\varepsilon,j})}\le K_p\varepsilon^{1/p}\big|\big|\phi_{\varepsilon}^j-i_{\partial C_{\varepsilon,j}}^*\bar F_{\varepsilon}^j\big|\big|_{L^p(\partial C_{\varepsilon,j})}\end{align}
holds for some positive constant $K_p>0$.\\
Eventually, we define the smooth $2$-form $\tilde F_{\varepsilon}^g$ on $C_{\varepsilon}^g$ by 
\[\tilde F_{\varepsilon}^g:= d\tilde A_{\varepsilon}^j \qquad \mbox{ on } C_{\varepsilon,j},\]
for every good cube $C_{\varepsilon,j}\subset C_{\varepsilon}^g$. We want to show that:
\begin{enumerate}
	\item $d\tilde F_{\varepsilon}^g=0$ on $\mathcal{D}'\big(\Int\big(\overline{C_{\varepsilon}^g}\big)\big)$. Without loss of generality, we are going to show that $d\tilde F_{\varepsilon}^g=0$ on $\mathcal{D}'\big(\Int\big(\overline{C_{\varepsilon,j}\cup C_{\varepsilon,k}}\big)\big)$, for every couple fo neighbouring good cubes $C_{\varepsilon,j}$ and $C_{\varepsilon,k}$. Indeed, pick any two neighbouring cubes $C_{\varepsilon,j},C_{\varepsilon,k}$ and define $\gamma_{jk}:=\partial C_{\varepsilon,j}\cap\partial C_{\varepsilon,k}$. Then, for every $\varphi\in C_c^{\infty}\big(\Int\big(\overline{C_{\varepsilon,j}\cup C_{\varepsilon,k}}\big)\big)$ it holds that
	\begin{align}0=\int_{C_{\varepsilon,j}\cup C_{\varepsilon,k}}d\big(\tilde F_{\varepsilon}^g\varphi\big)=\int_{C_{\varepsilon,j}\cup C_{\varepsilon,k}}d\tilde F_{\varepsilon}^g\varphi+\int_{C_{\varepsilon,j}\cup C_{\varepsilon,k}}\tilde F_{\varepsilon}^g\wedge d\varphi.\end{align}
	Thus, by Stokes theorem and (5.4) we deduce that
	\begin{align*}\int_{C_{\varepsilon,j}\cup C_{\varepsilon,k}}d\tilde F_{\varepsilon}^g\varphi&=-\int_{C_{\varepsilon,j}\cup C_{\varepsilon,k}}\tilde F_{\varepsilon}^g\wedge d\varphi=-\int_{C_{\varepsilon,j}}d\tilde A_{\varepsilon}^j\wedge d\varphi-\int_{C_{\varepsilon,k}}d\tilde A_{\varepsilon}^k\wedge d\varphi\\
	&=-\int_{C_{\varepsilon,j}}d\big(d\tilde A_{\varepsilon}^j\varphi\big)-\int_{C_{\varepsilon,k}}d\big(d\tilde A_{\varepsilon}^k\varphi\big)\\
	&=-\int_{\partial C_{\varepsilon,j}}i_{\partial C_{\varepsilon,j}}^*\big(d\tilde A_{\varepsilon}^j\big)\varphi-\int_{\partial C_{\varepsilon,k}}i_{\partial C_{\varepsilon,k}}^*\big(d\tilde A_{\varepsilon}^k\big)\varphi\\
	&=\pm\bigg(\int_{\gamma_{jk}}\phi_{\varepsilon}^j\varphi-\int_{\gamma_{jk}}\phi_{\varepsilon}^k\varphi\bigg)=0.
	\end{align*}
	\item It holds that
	\[\lim_{\varepsilon\rightarrow 0^+}\big|\big|\tilde F_{\varepsilon}^g-F\big|\big|_{L^p(C_{\varepsilon}^g)}=0.\]
	Indeed, by exploiting the elliptic estimate (5.3), we have that
	\begingroup
	\allowdisplaybreaks
	\begin{align*}\big|\big|\tilde F_{\varepsilon}^g-F\big|\big|_{L^p(C_{\varepsilon}^g)}^p&=\sum_{j=1}^{N_{\varepsilon}^g}\big|\big|d\tilde A_{\varepsilon}^j-F\big|\big|_{L^p(C_{\varepsilon,j})}^p\\
	&\le2^{p-1}\sum_{j=1}^{N_{\varepsilon}^g}\big|\big|d\tilde A_{\varepsilon}^j-\bar F_{\varepsilon}^j\big|\big|_{L^p(C_{\varepsilon,j})}^p+2^{p-1}\sum_{j=1}^{N_{\varepsilon}^g}\big|\big|\bar F_{\varepsilon}^j-F\big|\big|_{L^p(C_{\varepsilon,j})}^p\\
	&=2^{p-1}\sum_{j=1}^{N_{\varepsilon}^g}\big|\big|d\tilde A_{\varepsilon}^j-d\bar A_{\varepsilon}^j\big|\big|_{L^p(C_{\varepsilon,j})}^p+2^{p-1}\sum_{j=1}^{N_{\varepsilon}^g}\big|\big|\bar F_{\varepsilon}^j-F\big|\big|_{L^p(C_{\varepsilon,j})}^p\\
	&=2^{p-1}\sum_{j=1}^{N_{\varepsilon}^g}\big|\big|dA_{\varepsilon}^j\big|\big|_{L^p(C_{\varepsilon,j})}^p+2^{p-1}\sum_{j=1}^{N_{\varepsilon}^g}\big|\big|\bar F_{\varepsilon}^j-F\big|\big|_{L^p(C_{\varepsilon,j})}^p\\
	&\le2^{p-1}K_p\sum_{j=1}^{N_{\varepsilon}^g}\varepsilon\big|\big|\phi_{\varepsilon}^j-i_{\partial C_{\varepsilon,j}}^*\bar F_{\varepsilon}^j\big|\big|_{L^p(\partial C_{\varepsilon,j})}^p+2^{p-1}\sum_{j=1}^{N_{\varepsilon}^g}\big|\big|\bar F_{\varepsilon}^j-F\big|\big|_{L^p(C_{\varepsilon,j})}^p.\\
	\end{align*}
	\endgroup
	Then, we notice that
	\begin{align*}\sum_{j=1}^{N_{\varepsilon}^g}\varepsilon\big|\big|\phi_{\varepsilon}^j-i_{\partial C_{\varepsilon,j}}^*\bar F_{\varepsilon}^j\big|\big|_{L^p(\partial C_{\varepsilon,j})}^p&\le 2^{p-1}\Bigg(\sum_{j=1}^{N_{\varepsilon}^g}\varepsilon\big|\big|\phi_{\varepsilon}^j-i_{\partial C_{\varepsilon,j}}^* F\big|\big|_{L^p(\partial C_{\varepsilon,j})}^p+\sum_{j=1}^{N_{\varepsilon}^g}\varepsilon\big|\big|F-\bar F_{\varepsilon}^j\big|\big|_{L^p(\partial C_{\varepsilon,j})}^p\Bigg)\\
	&<2^{p-1}\Big(\varepsilon+\varepsilon\big|\big|F-\prescript{l}{}{\bar F}_{\varepsilon}\big|\big|_{L^p(\partial C_{\varepsilon})}^p+\varepsilon\big|\big|F-\prescript{r}{}{\bar F}_{\varepsilon}\big|\big|_{L^p(\partial C_{\varepsilon})}^p\Big).
	\end{align*}
	Thus, by Lemma \ref{L3.1}, it follows that
	\[\lim_{\varepsilon\rightarrow 0^+}\sum_{j=1}^{N_{\varepsilon}^g}\varepsilon\big|\big|\phi_{\varepsilon}^j-i_{\partial C_{\varepsilon,j}}^*\bar F_{\varepsilon}^j\big|\big|_{L^p(\partial C_{\varepsilon,j})}^p=0.\]
	On the other hand, by approximation through continuous functions and by absolute continuity, it follows that 
	\[\lim_{\varepsilon\rightarrow 0^+}\sum_{j=1}^{N_{\varepsilon}^g}\big|\big|\bar F_{\varepsilon}^j-F\big|\big|_{L^p(C_{\varepsilon,j})}^p=0.\]
\end{enumerate}
 \noindent
 Thus, we eventually have built a family of $2$-forms $\big\{\tilde F_{\varepsilon}^g\big\}_{\varepsilon\in R}$ on the connected open set $C_{\varepsilon}^g$ that are smooth on each open cube composing $C_{\varepsilon}^g$, continuous on all of $C_{\varepsilon}^g$ and such that 
 \begin{enumerate}
 	\item $d\tilde F_{\varepsilon}^g=0$ in $\mathcal{D}'\big(C_{\varepsilon}^g\big)$, for every $\varepsilon\in R$;
 	\item $\big|\big|\tilde F_{\varepsilon}^g-F\big|\big|_{L^p(C_{\varepsilon}^g)}\rightarrow 0$, as $\varepsilon\rightarrow 0^+$.
\end{enumerate}
\noindent

\section{Radial approximation on the bad cubes}
\noindent 
Let $C_{\varepsilon,j}$ be any bad cube in the cubic decomposition $C_{\varepsilon}$ of size $\varepsilon\in R$ and denote by $\pi_{\varepsilon,j}:C_{\varepsilon,j}\rightarrow\partial C_{\varepsilon,j}$ the radial projection for $C_{\varepsilon,j}$ to its boundary. Explicitly, the map $\pi_{\varepsilon,j}$ is given by 
\[\pi_{\varepsilon,j}(x):=\frac{\varepsilon}{2}\frac{x-x_{\varepsilon,j}}{||x-x_{\varepsilon,j}||_{*}}+x_{\varepsilon,j}, \qquad \mbox{ for every } x\in C_{\varepsilon,j},\]
where $x_{\varepsilon,j}$ denotes the centre of $C_{\varepsilon,j}$.\\
Moreover, we set
\[0\neq n_{\varepsilon,j}:=\frac{1}{2\pi}\int_{\partial C_{\varepsilon,j}}i_{\partial C_{\varepsilon,j}}^*F.\]
We define a $2$-form $F_{\varepsilon}^j$ on $C_{\varepsilon,j}$ by radial extension of the boundary datum $\phi_{\varepsilon}^j$ on $\partial C_{\varepsilon,j}$. Namely, we set
\[\tilde F_{\varepsilon}^j:=\pi_{\varepsilon,j}^*\phi_{\varepsilon}^j.\]
First, we claim that
\[d\tilde F_{\varepsilon}^j=2\pi n_{\varepsilon,j}\delta_{x_{\varepsilon,j}}dx^1\wedge dx^2\wedge dx^3, \qquad \mbox{ in } \mathcal{D}'\big(C_{\varepsilon,j}\big).\]
Indeed, we notice that for every $r\in(0,\varepsilon)$ the restriction of the map $\pi_{\varepsilon,j}$ to $\partial C_r(x_{\varepsilon,j})$ is a bi-Lipschitz map. For every $r\in(0,\varepsilon)$ we denote by $g_r:\partial C_{\varepsilon,j}\rightarrow \partial C_r(x_{\varepsilon,j})$ its inverse. Fixed any$\varphi\in C_c^{\infty}\big(C_{\varepsilon,j}\big)$, we define the function $h_y:(0,\varepsilon)\rightarrow\r$ by
\[h_y(r):=\varphi\big(g_r(y)\big), \qquad \mbox{ for every } r\in(0,\varepsilon),\]
for every $y\in\partial C_{\varepsilon,j}$. Then, we compute
\begin{align*}\left<\star d\tilde F_{\varepsilon}^j,\varphi\right>&=\left<\Div\big((\star \tilde F_{\varepsilon}^j)^{\sharp}\big),\varphi\right>=-\int_{C_{\varepsilon,j}}\big((\star\tilde F_{\varepsilon}^j)^{\sharp}\cdot\nabla\varphi\big)\, d\L^3\\
&=-\int_0^{\varepsilon}\bigg(\int_{\partial C_r(x_{\varepsilon,j})}\frac{1}{r^2}(\star\phi_{\varepsilon}^j)^{\sharp}\big(\pi_{\varepsilon,j}(x)\big)\cdot\nabla\varphi(x)\, d\H^2(x)\bigg)\, dr\\
&=-\int_0^{\varepsilon}\bigg(\int_{\partial C_{\varepsilon,j}}(\star\phi_{\varepsilon}^j)^{\sharp}(y)\cdot\nabla\varphi\big(g_r(y)\big)\, d\H^2(y)\bigg)dr\\
&=-\int_{\partial C_{\varepsilon,j}}(\star\phi_{\varepsilon}^j)^{\sharp}(y)\cdot\bigg(\int_0^{\varepsilon}\nabla\varphi\big(g_r(y)\big)\, dr\bigg)\, d\H^2(y)\\
&=-\int_{\partial C_{\varepsilon,j}}(\star\phi_{\varepsilon}^j)^{\sharp}(y)\cdot\nu_{\partial C_{\varepsilon,j}}(y)\bigg(\int_0^{\varepsilon}h_y'(r)\, dr\bigg)\, d\H^2(y)\\
&=\varphi(x_{\varepsilon,j})\int_{\partial C_{\varepsilon,j}}(\star\phi_{\varepsilon}^j)^{\sharp}(y)\cdot\nu_{\partial C_{\varepsilon,j}}(y)\, d\H^2(y)\\
&=\varphi(x_{\varepsilon,j})\int_{\partial C_{\varepsilon,j}}\phi_{\varepsilon}^j=2\pi n_{\varepsilon,j}\varphi(x_{\varepsilon,j})=\left<2\pi n_{\varepsilon,j}\delta_{x_{\varepsilon,j}},\varphi\right>.
\end{align*}
By arbitrariness of $\varphi\in C_c^{\infty}\big(C_{\varepsilon,j}\big)$, our claim follows.\\
We define a $2$-form in $\Omega_{p}^2\big(C_{\varepsilon}^b\big)$ by
\[\tilde F_{\varepsilon}^b:=\tilde F_{\varepsilon}^j \qquad \mbox{ on } C_{\varepsilon,j},\]
for every bad cube $C_{\varepsilon,j}\subset C_{\varepsilon}^b$.\\
By the same computation that we have performed in the previous section, we can show that
\[d\tilde F_{\varepsilon}^b=\Bigg(\sum_{j=1}^{N_{\varepsilon}^b}2\pi n_{\varepsilon,j}\delta_{x_{\varepsilon,j}}\Bigg)dx^1\wedge dx^2\wedge dx^3, \qquad \mbox{ in } \mathcal{D}'\big(C_{\varepsilon}^b\big).\]
Moreover, it holds that
\begin{align}\big|\big|\tilde F_{\varepsilon}^b-F\big|\big|_{L^p(C_{\varepsilon}^b)}\xrightarrow{\varepsilon\rightarrow 0^+} 0.\end{align}
Indeed, it holds that
\begin{align*}
\big|\big|\tilde F_{\varepsilon}^b-F\big|\big|_{L^p(C_{\varepsilon}^b)}^p&\le 2^{p-1}\Big(\big|\big|\tilde F_{\varepsilon}^b\big|\big|_{L^p(C_{\varepsilon}^b)}^p+\big|\big|F\big|\big|_{L^p(C_{\varepsilon}^b)}^p\Big)\\
&\le 2^{p-1}\Bigg(\varepsilon\big|\big|\phi_{\varepsilon}-i_{\partial C_{\varepsilon}}^*F\big|\big|_{L^p(\partial C_{\varepsilon})}^p+\varepsilon\sum_{j=1}^{N_{\varepsilon}^b}\big|\big|F\big|\big|_{L^p(\partial C_{\varepsilon,j})}^p+\big|\big|F\big|\big|_{L^p(C_{\varepsilon}^b)}^p\Bigg).
\end{align*}
By the Lebesgue point theorem and Fubini's theorem, for every bad cube $C_{\varepsilon,j}$ we can find a side-length $\varepsilon_j\in(0,\varepsilon)$ such that
\[\left|\big|\big|F\big|\big|_{L^p(\partial C_{\varepsilon,j})}^p-\big|\big|F\big|\big|_{L^p\big(\partial C_{\varepsilon_j}(x_{\varepsilon,j})\big)}^p\right|<\frac{1}{N_{\varepsilon}^b}\]
and
\[\big|\big|F\big|\big|_{L^p\big(\partial C_{\varepsilon_j}(x_{\varepsilon,j})\big)}^p\le\frac{1}{\varepsilon}||F||_{L^p(C_{\varepsilon,j})}^p.\]
In this way, by \ref{eq4.5} and Lemma \ref{L3.2} we get that
\begin{align*}
\big|\big|\tilde F_{\varepsilon}^b-F\big|\big|_{L^p(C_{\varepsilon}^b)}^p&\le 2^{p-1}\Big(\big|\big|\tilde F_{\varepsilon}^b\big|\big|_{L^p(C_{\varepsilon}^b)}^p+\big|\big|F\big|\big|_{L^p(C_{\varepsilon}^b)}^p\Big)\\
&\le 2^{p-1}\Bigg(\varepsilon\big|\big|\phi_{\varepsilon}-i_{\partial C_{\varepsilon}}^*F\big|\big|_{L^p(\partial C_{\varepsilon})}^p+\varepsilon+2\big|\big|F\big|\big|_{L^p(C_{\varepsilon}^b)}^p\Bigg)\xrightarrow{\varepsilon\rightarrow 0^+}0
\end{align*}
and (6.1) follows.
\noindent

\section{Final result and useful consequences}
\noindent
For every fixed $\varepsilon\in R$, define the set $\Gamma_{\varepsilon}$ given by $\Gamma_{\varepsilon}:=\partial\big(\overline{C_{\varepsilon}}\big)$. Clearly, there exists a Lipschitz deformation by retraction $H_{\varepsilon}:B_2\smallsetminus\overline{C_{\varepsilon}}\rightarrow\Gamma_{\varepsilon}$ from $B_2\smallsetminus\overline{C_{\varepsilon}}$ to $\Gamma_{\varepsilon}$. We define the following $2$-form 
\[\tilde F_{\varepsilon}:=\begin{cases}\tilde F_{\varepsilon}^g & \mbox{ on } C_{\varepsilon}^g\\\tilde F_{\varepsilon}^b & \mbox{ on } C_{\varepsilon}^b\\ \big(i_{\Gamma_{\varepsilon}}\circ H_{\varepsilon}\big)^*\phi_{\varepsilon} & \mbox{ on } B_2\smallsetminus\overline{C_{\varepsilon}},\end{cases}\]
where $i_{\Gamma_{\varepsilon}}:\Gamma_{\varepsilon}\rightarrow\partial C_{\varepsilon}$ is the inclusion map. 
We notice that, since $d\phi_{\varepsilon}=0$ in $L^{\infty}\big(\partial C_{\varepsilon}\big)$ and $\phi_{\varepsilon}$ is a Lipschitz $2$-form on $\partial C_{\varepsilon}$, then it holds that $\big(i_{\Gamma_{\varepsilon}}\circ H_{\varepsilon}\big)^*\phi_{\varepsilon}$ is a Lipschitz $2$-form on $B_2\smallsetminus\overline{C_{\varepsilon}}$ and 
\[d\Big(\big(i_{\Gamma_{\varepsilon}}\circ H_{\varepsilon}\big)^*\phi_{\varepsilon}\Big)=\big(i_{\Gamma_{\varepsilon}}\circ H_{\varepsilon}\big)^*(d\phi_{\varepsilon})=0 \qquad \mbox{ in } \mathcal{D}'\big(B_2\smallsetminus\overline{C_{\varepsilon}}\big).\]
By construction, as all the pieces are Lipschitz and coincide on the interfaces, it holds that $\tilde F_{\varepsilon}\in\Omega_{W^{1,\infty}}^2(B_2)$ and
\[d\tilde F_{\varepsilon}=\Bigg(\sum_{j=1}^{N_{\varepsilon}^b}2\pi n_{\varepsilon,j}\delta_{x_{\varepsilon,j}}\Bigg)dx^1\wedge dx^2\wedge dx^3 \qquad \mbox{ in } \mathcal{D}'(B_2).\] 
Moreover, by what we have showed so far and since $\L^3(B\smallsetminus\overline{C_{\varepsilon}})\rightarrow 0$ as $\varepsilon\rightarrow 0^+$, it holds that
\[\big|\big|\tilde F_{\varepsilon}-F\big|\big|_{L^p(B)}=\big|\big|\tilde F_{\varepsilon}^g-F\big|\big|_{L^p(C_{\varepsilon}^g)}+\big|\big|\tilde F_{\varepsilon}^b-F\big|\big|_{L^p(C_{\varepsilon}^b)}+\Big|\Big|\big(i_{\Gamma_{\varepsilon}}\circ H_{\varepsilon}\big)\phi_{\varepsilon}-F\Big|\Big|_{L^p(B\smallsetminus\overline{C_{\varepsilon}})}\xrightarrow{\varepsilon\rightarrow 0^+}0\]
Then, by standard convolution with a radial regularizing kernel, we can find a $2$-form $F_{\varepsilon}$ on $B$ such that
\[\big|\big|F_{\varepsilon}-\tilde F_{\varepsilon}\big|\big|_{L^p(B)}<\varepsilon,\]
\[F_{\varepsilon}\in C^{\infty}\Bigg(B\smallsetminus\bigcup_{j=1}^{N_{\varepsilon}^b}\{x_{\varepsilon,j}\}\Bigg)\]
and
\[dF_{\varepsilon}=\Bigg(\sum_{j=1}^{N_{\varepsilon}^b}2\pi n_{\varepsilon,j}\delta_{x_{\varepsilon,j}}\Bigg)dx^1\wedge dx^2\wedge dx^3 \qquad \mbox{ in } \mathcal{D}'(B).\] 
This of course ensures also that
\[\big|\big|F_{\varepsilon}-F\big|\big|_{L^p(B)}\xrightarrow{\varepsilon\rightarrow 0^+}0\]
and the statement of Theorem \ref{T1.1} follows. 

\subsection{Characterization of the class $L_{\z}^p(B)$ when $p\ge 3/2$}
\noindent \smallskip \\
The proof of Corollary \ref{C1.2} is a direct consequence of the strong approximation result for elements of $\Omega_{p,\z}^2(B)$, Remark \ref{R1.3} and the following lemma. 
\begin{lem}
\label{L7.1}
Let $n>1$, $p\ge n/(n-1)$ and $X$ be any vector field in $L^p(B,\r^n)$, where $B\subset\r^n$ is the open unit ball in $\r^n$.\\
If the distributional divergence of $X$ is a finite linear combination of delta distributions, then $X$ is divergence free.\\
In particular
\[\Omega_{p,R}^2(B)\subset\mathscr{Z}_p^2(B),\]
where $\mathscr{Z}_p^2(B)$ is the set of the $L^p$ 2-forms $F$ on $B$ such that $dF=0$ in $\mathcal{D}'(B)$.
\begin{proof}
Without loss of generality, we assume that
\[\Div(X)=\alpha\delta_0 \qquad \mbox{ in } \mathcal{D}'(B) \mbox{ for some } \alpha\in\r\]
and $p=n/(n-1)$, which in turn implies $p'=n$.\\
Consider the sequence $\{\varphi_k\}_{k\in\n}$ of Lipschitz and compactly supported functions on $B$ defined by 
\[\varphi_k(x):=\begin{cases}k & \mbox{ for every } x\in\ B_{e^{-k}/2},\\-\ln\big(2|x|\big) & \mbox{ for every } x\in B_{1/2}\smallsetminus\overline{B_{e^{-k}/2}},\\ 0 & \mbox{ for every } x\in B\smallsetminus\overline{B_{1/2}}.\end{cases}\]
It follows easily that
\[|\nabla\varphi_k(x)|:=\begin{cases}0 & \mbox{ for every } x\in\ B_{e^{-k}/2},\\\displaystyle{\frac{1}{|x|}} & \mbox{ for every } \  x\in B_{1/2}\smallsetminus B_{e^{-k}/2},\\ 0 & \mbox{ for every } x\in B\smallsetminus B_{1/2}.\end{cases}\]
We then compute
\begin{align*}||\nabla\varphi_k||_{L^n(B)}&=\bigg(\int_{B}|\nabla\varphi_k(x)|^n\, d\L^3(x)=\bigg)^{1/n}=\bigg(\int_{B_{1/2}\smallsetminus B_{e^{-k}/2}}|x|^{-n}\, d\L^3(x)=\bigg)^{1/n}\\
&=\bigg(\int_{e^{-k}/2}^{1/2}r^{-n}\int_{\partial B_{r}}d\H^{n-1}(y)\, dr\bigg)^{1/n}=\bigg(\int_{e^{-k}/2}^{1/2}r^{-n}n\omega_nr^{n-1}\, dr\bigg)^{1/n}\\
&=(n\omega_n)^{1/n}\bigg(\int_{e^{-k}/2}^{1/2}r^{-1}\bigg)^{1/n}=(n\omega_n)^{1/n}\Big(\big[\ln(r)\big]_{e^{-k}/2}^{1/2}\Big)^{1/n}\\
&=(n\omega_n)^{1/n}\big(-\ln(2)+\ln(2)-\ln(e^{-k})\big)^{1/n}=(n\omega_n)^{1/n}k^{1/n}
\end{align*}
and we notice that
\begin{align*}|\alpha|k&=\left|\left<\Div(X),\varphi_k\right>\right|\le\int_{B}(|X|\cdot|\nabla\varphi_k|\, d\L^3\le||X||_{L^{n/(n-1)}(B)}||\nabla\varphi_k||_{L^n(B)}\\
&=\Big((n\omega_n)^{1/n}||X||_{L^{n/(n-1)}(B)}\Big)k^{1/n}=C_n||X||_{L^{n/(n-1)}(B)}k^{1/n},
\end{align*}
where $C_n:=(n\omega_n)^{1/n}$. By letting $k\rightarrow+\infty$ we get $|\alpha|=0$ and the statement follows.
\end{proof}
\end{lem}
\noindent 

\subsection{A decomposition theorem for vector fields in $L_{\z}^1(B)$}
\noindent \smallskip \\
This subsection is entirely dedicated to the proof of Theorem \ref{T1.2}. Such theorem states that, given any vector field $X\in L_{\z}^1(B)$, we can identify $X$ with a $1$-current on $B$ which can be written as the sum of $1$-cycle and a calibrated integer multiplicity rectifiable $1$-current.\\
First of all, notice that we can associate to any $X\in L^1(B;\r^3)$ a $1$-current on $B$ given 
\[\left<T_X,\omega\right>:=\int_{B}\left<\omega,X\right>\, d\L^3, \qquad \mbox{ for every } \omega\in\mathcal{D}^1(B).\]
Let's start to study the structure of $T_X$ for elements in $L_R^1(B)$, i.e. for smooth vector fields up to finitely many integer degree singularities on $B$.
\begin{lem}
\label{L7.2}
Let $X\in L_R^1(B)$. Then, there exists an integer $1$-current $L\in\mathcal{R}_1(B)$ with finite mass and a $1$-cycle $C\in\mathcal{D}_1(B)$ such that $T_X=C+L$.
\begin{proof}
We claim that there exists an integer multiplicity $1$-current $L\in\mathcal{R}_1(B)$ such that $\partial T_X=\partial L$. Then, the statement will follows by setting $C:=T_X-L$.\\
In order to prove our claim, we proceed as follows. By definition, it holds that $\partial T_X=\Div(X)$ in $\mathcal{D}'(B)$ and 
\[\Div(X)=\sum_{j=1}^nd_j\delta_{x_j}, \qquad \mbox{ for some } d_1,...,d_n\in\z\smallsetminus\{0\} \mbox{ and } x_1,...,x_n\in B.\] 
We let $p,q\in\n$ be the number of points $x_j$ in $\{x_1,...,x_n\}$ such that $d_j$ is positive and negative respectively. Then, we define 
\[\{i_1,...,i_p\}=\big\{j\in\{1,...,n\} \mbox{ s.t. } d_j>0\big\},\]
\[\{j_1,...,j_q\}=\big\{j\in\{1,...,n\} \mbox{ s.t. } d_j<0\big\}\]
and
\[d=\sum_{j=1}^nd_j\in\z.\]
We build a family $\mathscr{F}=\{L_{\alpha}\}_{\alpha\in A}$ of oriented segments in $B$ as follows. If $q=0$, we set $\mathscr{F}=\emptyset$. If $q>0$, we start from $x_{i_1}$ and we add to the family $\mathscr{F}$ the segments $\big(x_{j_1},x_{i_1}\big),...,\big(x_{j_{k_1}},x_{i_1}\big)$, until we reach the condition $k_1=q$ or the condition 
\[r_1:=d_{i_1}+\sum_{l=1}^{k_1}d_{j_l}\le 0.\]
If $k_1=q$, then we stop. If $k_1<q$, then we move to $x_{i_2}$. If $r_1=0$, then we add to $\mathscr{F}$ the segments $\big(x_{j_{k_1+1}},x_{i_2}\big),...,\big(x_{j_{k_2}},x_{i_2}\big)$, where $k_2\in\{1,...,q\}$ is the smallest value such that
\[r_2:=d_{i_2}+\sum_{l=k_1+1}^{k_2}d_{j_l}\le 0.\]
If $r_1<0$, then we add to $\mathscr{F}$ the segment $\big(x_{j_{k_1}},x_{i_2}\big)$ and the segments $\big(x_{j_{k_1+1}},x_{i_2}\big),...,\big(x_{j_{k_2}},x_{i_2}\big)$, where $k_2\in\{1,...,q\}$ is the smallest value such that
\[r_2:=d_{i_2}+r_1+\sum_{l=k_1+1}^{k_2}d_{j_l}\le 0.\]
We proceed iteratively in this way, moving on the subsequent nodes $x_{i_s}$ until the counter $k_s=q$ or $s=p$. Then, the construction of the family $\mathscr{F}$ is complete. We let $x_{i_h}$ be the last node $x_{i_s}$ that is visited before the iteration stops and, for every $L_{\alpha}=(x_j,x_i)\in\mathscr{F}$, we define its multiplicity $m_{\alpha}$ as 
\[m_{\alpha}:=\min\{|d_j|,d_i|\}.\] 
Eventually, we divide three cases:
\begin{enumerate}
\item \textit{Case} $d=0$. We define the integer $1$-current 
$L \in\mathcal{R}_1(B)$ given by
\[\left<L,\omega\right>:=\sum_{\alpha\in A}m_{\alpha}\int_{L_{\alpha}}\omega, \qquad \mbox{ for every } \omega\in\mathcal{D}^1(B).\]
\item \textit{Case} $d>0$. We fix a point $x_0\in\partial B$ and we let $L_s^p:=(x_0,x_{i_s})$, for every $s=h,...,p$. We define the integer $1$-current 
$L \in\mathcal{R}_1(B)$ given by
\[\left<L,\omega\right>:=\sum_{\alpha\in A}m_{\alpha}\int_{L_{\alpha}}\omega+r_h\int_{L_h^b}\omega+\sum_{s=h+1}^pd_{i_s}\int_{L_s^b}\omega,\qquad \mbox{ for every } \omega\in\mathcal{D}^1(B).\]
\item \textit{Case} $d<0$. We fix a point $x_0\in\partial B$ and we let $L_s^n:=(x_{j_s},x_0)$, for every $s=k_h,...,q$ We define the integer $1$-current $L \in\mathcal{R}_1(B)$ given by
\[\left<L,\omega\right>:=\sum_{\alpha\in A}m_{\alpha}\int_{L_{\alpha}}\omega+|r_h|\int_{L_h^b}\omega+\sum_{s=k_h+1}^q|d_{i_s}|\int_{L_s^b}\omega,\qquad \mbox{ for every } \omega\in\mathcal{D}^1(B).\]
\end{enumerate}
By direct computation it can be shown that $L$ has the desired properties. Hence, the statement follows.
\end{proof}
\end{lem}
We
\begin{lem}
\label{L7.3}
Let $X\in L_R^1(B)$. Then, 
\begin{align}
\inf_{\substack{T\in\mathcal{D}_1(B), \\ \partial T=\partial T_X}}\mathbb{M}(T)=\inf_{\substack{T\in\mathcal{M}_1(B), \\ \partial T=\partial T_X}}\mathbb{M}(T)=\sup_{\substack{\varphi\in\mathcal{D}(B), \\ ||\nabla\varphi||_{L^{\infty}(B)}\le 1}}\left<\partial T_X,\varphi\right>,
\end{align}
where $\mathcal{M}_1(B)$ denotes the set of all the $1$-currents with finite mass on $B$.
\begin{proof}
First of all, by Lemma \ref{L7.2} we infer that there exists an integer $1$-current $L\in\mathcal{R}_1(B)$ with finite mass such that $\partial L=\partial T_X$. Hence, 
\[\inf_{\substack{T\in\mathcal{D}_1(B), \\ \partial T=\partial T_X}}\mathbb{M}(T)\le M(L)<+\infty\]
and we conclude that
\[\inf_{\substack{T\in\mathcal{D}_1(B), \\ \partial T=\partial T_X}}\mathbb{M}(T)=\inf_{\substack{T\in\mathcal{M}_1(B), \\ \partial T=\partial T_X}}\mathbb{M}(T).\]
Since for every $T\in\mathcal{M}_1(B)$ such that $\partial T=\partial T_X$ it holds 
\[\left<\partial T_X,\varphi\right>=\left<\partial T,\varphi\right>=\left<T,d\varphi\right>\le\mathbb{M}(T)||\nabla\varphi||_{L^{\infty}(B)}, \qquad \mbox{ for every } \varphi\in\mathcal{D}(B),\]
then 
\begin{align}\inf_{\substack{T\in\mathcal{M}_1(B), \\ \partial T=\partial T_X}}\mathbb{M}(T)\ge\sup_{\substack{\varphi\in\mathcal{D}(B), \\ ||\nabla\varphi||_{L^{\infty}(B)}\le 1}}\left<\partial T_X,\varphi\right>.
\end{align}
To prove that the inequality (7.2) is actually an equality, it's sufficient to show that the supremum on the right-hand side of (7.2) equals the mass of some $1$-current with finite mass $T$. To build such a current, we first consider the following linear subspace of $C^0(\overline{B};\r^3)$:
\[D:=\big\{\xi\in C^0(\overline{B};\r^3) \mbox{ such that } \xi|_B=\nabla\varphi \mbox{ for some } \varphi\in\mathcal{D}(B)\big\}.\]
Then, we define the linear functional $f:D\subset C^0(\overline{B};\r^3)\rightarrow\r$ by
\[\left<f,\xi\right>=\left<\partial T_X,\varphi\right>, \mbox{ where } \xi=\nabla\varphi, \mbox{ for every } \xi\in D.\]
Since
\[\left<f,\xi\right>=\left<\partial T_X,\varphi\right>=\left<\partial L,\varphi\right>=\left<L,d\varphi\right>\le\mathbb{M}(L)||\nabla\varphi||_{L^{\infty}(B)}=\mathbb{M}(L)||\xi||_{L^{\infty}(B)}, \qquad \mbox{ for every } \xi\in D,\]
then $f\in D^*$. By the Hahn-Banach theorem, we can extend $f$ to a continuous linear functional $F$ on $C^0(\overline{B};\r^3)$ by preserving its norm. We define the$1$-current $T:\mathcal{D}^1(B)\rightarrow\r$ by
\[\left<T,\omega\right>=\left<F,\xi^{\omega}\right>,\]
where $\xi^{\omega}\in C^0(B)$ is the extension of the vector field $\omega^{\sharp}$ by zero on the boundary of $B$, for every $\omega\in\mathcal{D}^1(B)$.\\
By definition, it holds that $\partial T=\partial T_X$ and, moreover, it holds that
\[\mathbb{M}(T)=||F||_{C^0(\overline{B},\r^3)^*}=||f||_{D*}=\sup_{\substack{\xi\in D, \\ ||\xi||_{L^{\infty}(B)}\le 1}}\left<f,\xi\right>=\sup_{\substack{\varphi\in\mathcal{D}(B), \\ ||\nabla\varphi||_{L^{\infty}(B)}\le 1}}\left<\partial T_X,\varphi\right><+\infty.\]
Thus, the statement follows and we also conclude that the infimum on the left-hand side is achieved by some $1$-current with finite mass on $B$. 
\end{proof}
\end{lem}
\begin{thm}
\label{T7.1}
Let $X\in L_R^1(B)$. Then, there exists an integer $1$-current $L\in\mathcal{R}_1(B)$ with finite mass and a $1$-cycle $C\in\mathcal{D}_1(B)$ such that:
\begin{enumerate}
\item $T_X=C+L$;
\item $L$ is a mass-minimizer on the class $\mathcal{R}_1(B)\cap\{T\in\mathcal{D}_1(B) \mbox{ s.t. } \partial T=\partial T_X\}$ and its mass is given by
\[\mathbb{M}(L)=\sup_{\substack{\varphi\in\mathcal{D}(B), \\ ||\nabla\varphi||_{L^{\infty}(B)}\le 1}}\left<\partial T_X,\varphi\right>;\]
\item $L$ is calibrated by an exact, essentially bounded $1$-form $\omega=d\varphi$, with $\varphi\in\Lip_1(B)$.
\end{enumerate}
\begin{proof}
First of all, notice that by Lemma \ref{L7.2}, Lemma \ref{L7.3} and \cite[Charper 1, Section 3.4, Theorem 8]{cartesiancurrents2} it follows that
\[\inf_{\substack{T\in\mathcal{R}_1(B), \\ \partial T=\partial T_X}}\mathbb{M}(T)=\sup_{\substack{\varphi\in\mathcal{D}(B), \\ ||\nabla\varphi||_{L^{\infty}(B)}\le 1}}\left<\partial T_X,\varphi\right><+\infty.\]
Since the mass $\mathbb{M}(\cdot)$ is lower semicontinuous with respect to the weak topology on $\mathcal{D}_1(B)$ and the competition class $\mathcal{R}_1(B)\cap\{T\in\mathcal{D}_1(B) \mbox{ s.t. } \partial T=\partial T_X\}$ is compact with respect to the same topology (for a reference, see e.g. \cite[Equation (7.5), Theorem 7.5.2]{krantzparks}), by the direct method of calculus of variations we conclude that there exists an integer $1$-current $L\in\mathcal{R}_1(B)$ such that $\partial L=\partial T_X$ and 
\[\mathbb{M}(L)=\inf_{\substack{T\in\mathcal{R}_1(B), \\ \partial T=\partial T_X}}\mathbb{M}(T)=\sup_{\substack{\varphi\in\mathcal{D}(B), \\ ||\nabla\varphi||_{L^{\infty}(B)}\le 1}}\left<\partial T_X,\varphi\right>.\]
This is enough to prove (1) and (2). For what concerns (3), we consider a sequence $\{\varphi_k\}_{k\in\n}\subset\mathcal{D}(B)$ such that $||\nabla\varphi_k||_{L^{\infty}(B)}\le 1$ for every $k\in\n$ and 
\[\lim_{k\rightarrow+\infty}\left<\partial T_X,\varphi_k\right>=\mathbb{M}(L).\]
Since $\{\varphi_k\}_{k\in\n}$ is a sequence of uniformly $1$-Lipschitz functions on $\overline{B}$, by the Ascoli-Arzelà theorem we can find a subsequence $\big\{\varphi_{k_h}\big\}_{h\in\n}$ such that $\varphi_{k_h}\xrightarrow{h\rightarrow\infty} \varphi$ uniformly on $B$, for some $1$-Lipschitz function $\varphi\in\Lip_1(B)$. We set $\omega:=d\varphi$ and we want to show that
\[\left<L,\omega\right>=\mathbb{M}(L).\]
Indeed,
\[\left<L,\omega\right>=\left<L,d\varphi\right>=\left<\partial L,\varphi\right>=\left<\partial T_X,\varphi\right>=\lim_{k\rightarrow+\infty}\left<\partial T_X,\varphi_k\right>=\mathbb{M}(L).\]
Thus, the statement follows.
\end{proof}
\end{thm}
\noindent
Now we are ready to prove Theorem \ref{T1.2}, basing on the previous lemmata and on the strong approximation result for vector fields in $L_{\z}^1(B)$.
\begin{proof}[Proof of Theorem 1.2]
Let $X\in L_{\z}^1(B)$. By the strong approximation theorem, we can pick a sequence of vector fields $\{X_k\}_{k\in\n}\subset L_R^1(B)$ such that $X_k\xrightarrow{k\rightarrow\infty} X$ in $L^1(B)$. Since the estimate
\[\left|\left<\big(T_{X_k}-T_{X}\big),\omega\right>\right|=\int_{B}\left<\omega,X_k-X\right>, d\L^3\le ||\omega||_{L^{\infty}(B)}||X_k-X||_{L^1(B)}\le ||X_k-X||_{L^1(B)}\]
holds for every $\omega\in\mathcal{D}^1(B)$ such that $||\omega||_{L^{\infty}(B)}\le 1$, we conclude that
\[\mathbb{M}\big(T_{X_k}-T_X\big)=\sup_{\substack{\omega\in\mathcal{D}^1(B), \\ ||\omega||_{L^{\infty}(B)}\le 1}}\left<T_{X_k}-T_X,\omega\right>\le||X_k-X||_{L^1(B)}\xrightarrow{k\rightarrow\infty} 0.\]
Moreover, since
\[\sup_{\substack{\varphi\in\mathcal{D}(B), \\ ||\nabla\varphi||_{L^{\infty}(B)}\le 1}}\left<\partial T_{X_k}-\partial T_X,\varphi\right>=\sup_{\substack{\varphi\in\mathcal{D}(B), \\ ||\nabla\varphi||_{L^{\infty}(B)}\le 1}}\left<T_{X_k}-T_X,d\varphi\right>\le\mathbb{M}\big(T_{X_k}-T_X\big),\]
we deduce that
\[\sup_{\substack{\varphi\in\mathcal{D}(B), \\ ||\nabla\varphi||_{L^{\infty}(B)}\le 1}}\left<\partial T_{X_k}-\partial T_X,\varphi\right>\xrightarrow{k\rightarrow\infty}0.\]
Thus, fixed any $0<\varepsilon<1$, we can find a subsequence $\big\{X_{k_j(\varepsilon)}\big\}_{j\in\n}\subset L_R^1(B)$ such that
\[\sup_{\substack{\varphi\in\mathcal{D}(B), \\ ||\nabla\varphi||_{L^{\infty}(B)}\le 1}}\left<\partial T_{X_{k_j(\varepsilon)}}-\partial T_{X_{k_{j+1}(\varepsilon)}},\varphi\right>\le\frac{\varepsilon}{2^j}, \qquad \mbox{ for every } j\in\n.\]
For every $j\in\n$, let $L_j^{\varepsilon}$ be the minimal integer $1$-current such that $\partial L_j^{\varepsilon}=\partial T_{X_{k_j(\varepsilon)}}$ whose existence is given by Theorem \ref{T7.1}. Analogously, for every $j\in\n$, let $L_{j,j+1}^{\varepsilon}$ be the the minimal integer $1$-current such that $\partial L_{j,j+1}^{\varepsilon}=\partial T_{X_{k_j(\varepsilon)}}-\partial T_{X_{k_{j+1}(\varepsilon)}}=\partial L_j^{\varepsilon}-\partial L_{j+1}^{\varepsilon}$, who existence is again given by Theorem \ref{L7.1}.\\
Define the following sequence of integer $1$-currents on $B$:
\[L_n^{\varepsilon}:=\begin{cases}
L_0^{\varepsilon} & \mbox{ if } n=0,\\ \displaystyle{L_0^{\varepsilon}-\sum_{j=0}^{n-1}L_{k_jk_{j+1}}^{\varepsilon}}& \mbox{ if } n>0,
\end{cases} \qquad \mbox{ for every } n\in\n.\]
Clearly, 
\begin{align}\partial L_n^{\varepsilon}=\partial L_0^{\varepsilon}-\sum_{j=0}^{n-1}\partial L_{j,j+1}^{\varepsilon}=\partial L_0^{\varepsilon}-\sum_{j=0}^{n-1}(\partial L_j^{\varepsilon}-\partial L_{j+1}^{\varepsilon})=\partial L_n^{\varepsilon}=\partial T_{X_{k_n(\varepsilon)}}.
\end{align}
Moreover, since $L_{j,j+1}^{\varepsilon}$ is minimal for every $j\in\n$, it holds that
\[\mathbb{M}\big(L_{j,j+1}^{\varepsilon}\big)=\sup_{\substack{\varphi\in\mathcal{D}(B), \\ ||\nabla\varphi||_{L^{\infty}(B)}\le 1}}\left<\partial T_{X_{k_j(\varepsilon)}}-\partial T_{X_{k_{j+1}(\varepsilon)}},\varphi\right>\le\frac{\varepsilon}{2^j}, \qquad \mbox{ for every } j\in\n.\]
Thus, the estimate
\[\mathbb{M}(L_{n+1}^{\varepsilon}-L_n^{\varepsilon})=\mathbb{M}\big(L_{n,n+1}^{\varepsilon}\big)\le\frac{\varepsilon}{2^n}, \qquad \mbox{ for every } n\in\n,\]
let us conclude that the sequence $\{L_n^{\varepsilon}\}_{n\in\n}$ is a Cauchy sequence in mass. Hence, there exists an integer $1$-current $L^{\varepsilon}\in\mathcal{R}_1(B)$ such that 
\begin{align}\mathbb{M}(L_n^{\varepsilon}-L^{\varepsilon})\xrightarrow{n\rightarrow\infty}0,\end{align}
By (7.3), the strong convergence (7.4) and $X_k\xrightarrow{k\rightarrow\infty}X$ in $L^1(B)$, we conclude that
\[\partial L^{\varepsilon}=\lim_{n\rightarrow +\infty}\partial L_n^{\varepsilon}=\lim_{n\rightarrow+\infty}\partial T_{X_{k_n(\varepsilon)}}=\partial T_X.\]
Since the family of integer $1$-cycles $\{L^{\varepsilon}-L^{1/2}\}_{0<\varepsilon<1}\subset\mathcal{R}_1(B)$ is uniformly bounded in mass, then by standard compactness arguments for currents we can find a sequence $\varepsilon_k\rightarrow 0$ as $k\rightarrow+\infty$ and an integer 1-cycle $\tilde L\in\mathcal{R}_1(B)$ with finite mass such that $L^{\varepsilon_k}-L^{1/2}\rightharpoonup \tilde L$. If we let $L:=L^{1/2}+\tilde L$ we get $L^{\varepsilon_k}\rightharpoonup L$ weakly in $\mathcal{D}_1(B)$. By construction, $L$ is again an integer $1$-current with finite mass such that $\partial L=\partial L^{1/2}=\partial T_X$. We claim that
\[\mathbb{M}(L)=\inf_{\substack{T\in\mathcal{M}_1(B), \\ \partial T=\partial T_X}}\mathbb{M}(T).\] 
By contradiction, assume that we can find $T\in\mathcal{M}_1(B)$ such that $\partial T=\partial T_X$ and 
\begin{align}\mathbb{M}(T)<\mathbb{M}(L)\le\liminf_{k\rightarrow+\infty}\mathbb{M}\big(L^{\varepsilon_k}\big),\end{align}
where the last inequality follows by weak convergence and lower semicontinuity of the mass.\\
By (7.5) we can find some $h\in\n$ such that
\begin{align}\mathbb{M}(T)<\mathbb{M}\big(L^{\varepsilon_{h}}\big)-2\varepsilon_{h}.\end{align}
Moreover, since $\mathbb{M}\big(L_0^{\varepsilon}-L^{\varepsilon}\big)\le\varepsilon$ for every $0<\varepsilon<1$, it holds that
\begin{align}\mathbb{M}\big(L_0^{\varepsilon_h}-L^{\varepsilon_h}\big)\le\varepsilon_{h}.\end{align}
We define $\tilde T:=T+L_0^{\varepsilon_h}-L^{\varepsilon_h}$ and we notice that $\partial\tilde T=\partial T_{X_{k_0}}$. Moreover, by (7.6), (7.7) and the minimality of $L_0^{\varepsilon_h}$, we conclude that
\[\mathbb{M}\big(L_0^{\varepsilon_h}\big)\le\mathbb{M}\big(\tilde T\big)\le\mathbb{M}(T)+\mathbb{M}\big(L_0^{\varepsilon_h}-L^{\varepsilon_h}\big)<\mathbb{M}\big(L^{\varepsilon_h}\big)-\varepsilon_h\le\mathbb{M}\big(L_0^{\varepsilon_h}\big),\]
which is a contradiction. Thus, our claim follows.\\
Since $L\in\mathcal{R}_1(B)$, we get that
\[\mathbb{M}(L)=\inf_{\substack{T\in\mathcal{R}_1(B), \\ \partial T=\partial T_X}}\mathbb{M}(T)=\inf_{\substack{T\in\mathcal{M}_1(B), \\ \partial T=\partial T_X}}\mathbb{M}(T)\]
and, by repeating the same proof as in Lemma \ref{L7.3}, we obtain that
\[\inf_{\substack{T\in\mathcal{M}_1(B), \\ \partial T=\partial T_X}}\mathbb{M}(T)=\sup_{\substack{\varphi\in\mathcal{D}(B), \\ ||\nabla\varphi||_{L^{\infty}(B)}\le 1}}\left<\partial T_X,\varphi\right>.\]
Hence, by setting $C:=T_X-L$ the statement follows immediately. 
\end{proof}
\section*{Acknowledgements}
\noindent
I want to include this section in order to thank my Ph.D. advisor, Prof. Dr. Tristan Rivière, for having encouraged me to think about these problems and for the mathematical and personal support that he provided me when my research was stuck during this first year at the ETH Zürich.\\
Moreover, I'd like to say thank you to my colleagues and friends Federico Franceschini, Giada Franz, Federico Glaudo and Dr. Alessandro Pigati, for our useful discussions and for the constant feedback that they've been giving me about my job. 
\appendix
\section{Harmonic extensions of $1$-forms on open cubes}
\begin{lem}
\label{LA.1}
Let $Q\subset\r^n$ be an open cube in $\r^n$ with centre $x_0$ and side-length $r>0$. Choose any $f\in C_c^{\infty}(\r^n)$ and let $u\in W^{1,2}(Q)$ be a weak solution of the following differential problem:
\begin{align}\begin{cases}\Delta u=0&\mbox{ on } Q\\u=f|_{\partial Q} &\mbox{ on } \partial Q.\end{cases}\end{align}
Then, for every $q\in [1,+\infty)$ it holds that $u\in W^{2,q}(Q)\cap C^{\infty}(Q)$. Moreover, for every $q\in(1,+\infty)$ there exists a constant $K_q^0>0$ such that
\[||\nabla u||_{L^q(Q)}\le K_q^0||\nabla f||_{L^q(Q)}.\]
\begin{proof}
Consider the function $v:=u+f$. Then, $v\in W_0^{1,2}(Q)$ is a weak solution of the differential problem
\begin{align}\begin{cases}\Delta v=\Delta f=:h&\mbox{ on } Q\\v=0 &\mbox{ on } \partial Q.\end{cases}\end{align}
Since $h\in C^{\infty}(Q)$, by standard elliptic interior regularity we immediately conclude that $v\in C^{\infty}(Q)$. In order to recover at lest a partial regularity for $u$ up to the boundary of $Q$, we argue as follows.\\ 
Given any function $g:Q\rightarrow\r$, we define an extension $\tilde g$ of $g$ to the open cube $Q'=\overline{Q}+(-l,l)^n\supset Q$ by simply reflecting $g$ with respect to the faces, the sides and the vertices of $\partial Q$. In such a way, we get that the function $\tilde v\in W_0^{1,2}(Q)$ solves the differential problem
\[\begin{cases}\Delta\tilde v=\tilde h&\mbox{ on } Q'\\\tilde v=0 &\mbox{ on } \partial Q'.\end{cases}\]
By standard elliptic regularity, it follows that $\tilde v|_{\Omega}\in W^{2,q}(\Omega)$, for every $q\in[1,+\infty)$ for every $\Omega\subset\subset Q'$. Hence, it follows that $v=\tilde v|_Q\in W_0^{2,q}(Q)\cap C^{\infty}(Q)$, for very $q\in[1,+\infty)$. Since $f\in C^{\infty}(\overline{Q})\subset W^{2,p}(Q)\cap C^{\infty}(Q)$, we eventually conclude that $u=v-f\in W^{2,q}(Q)\cap C^{\infty}(Q)$, for every $q\in[1,+\infty)$. In particular, by Sobolev embedding theorem, it follows that $u\in C^{1,\alpha}(\overline{Q})$, for every $0\le\alpha<1$.\\
To conclude the proof, we just have to show the estimate (A.1). By applying the result in section 3.5 of \cite{amb-car-mas} to the differential problem (A.2) we obtain that for every $q\in(1,+\infty)$ there exists a constant $\tilde K(p,Q)>0$ such that
\[||\nabla v||_{L^p(Q)}\le\tilde K(q,Q)||\nabla f||_{L^p(Q)}.\]
Thus, by the triangular inequality, we get that
\[||\nabla u||_{L^q(Q)}\le K(q,Q)||\nabla f||_{L^q(Q)}.\]
with $K(q,Q):=1+\tilde K(q,Q)$. Let $K_q^0:=K(q,C)$, where $C$ is the open unit cube in $\r^n$ and observe that the function $\tilde u:C\rightarrow\r$ given by $\tilde u:=u\big(r(\hspace{0.5mm}\cdot-x_0)\big)$ is a weak solution in $W^{2,q}(C)\cap C^{\infty}(C)$ for the differential problem
\begin{align*}\begin{cases}\Delta\tilde u=0&\mbox{ on } C\\\tilde u=\tilde f|_{\partial C} &\mbox{ on } \partial C,\end{cases}\end{align*}
where $\tilde f=f\big(r(\hspace{0.5mm}\cdot-x_0)\big)\in C_c^{\infty}(\r^n)$. Then, by what we have proved so far, it holds that
\[||\nabla u||_{L^q(Q)}=r^{-\left(1-\frac{n}{q}\right)}||\nabla\tilde u||_{L^q(C)}\le r^{-\left(1-\frac{n}{q}\right)}K_q^0||\nabla\tilde f||_{L^q(C)}=K_q^0||\nabla f||_{L^q(Q)}\]
and the statement follows.
\end{proof}
\end{lem}
\begin{lem}
\label{LA.2}
Let $Q\subset\r^n$ be an open cube in $\r^n$ with centre $x_0$ and side-length $r>0$. Let $B\in\Omega_c^1(\r^n)$, $f\in L^{\infty}(\partial Q)$ and assume that the $1$-form $\beta:=i_{\partial Q}^*B\in\Omega_{W^{1,\infty}}^1(Q)$ is such that
\begin{align}\begin{cases}d\beta=f\\d^*\beta=0\end{cases}\mbox{ on } \partial Q.\end{align}
Then:
\begin{enumerate}
\item there exists an exponent $\tilde q\in(1,+\infty)$ and a constant $\kappa_1(Q)>0$ such that 
\[||B||_{W^{1,\tilde q}(Q)}\le\kappa_1(Q)||f||_{L^1(\partial Q)};\]
\item for every $q\in(1,+\infty)$ there exists a constant $\kappa_q(Q)>0$ such that
\[||B||_{W^{1,q}(Q)}\le\kappa_q(Q)||f||_{L^q(\partial Q)}.\]
\end{enumerate}
\begin{proof}
As the trace operator $i_{\partial Q}^*:W^{1,q}(Q)\rightarrow W^{1-\frac{1}{q},q}(\partial Q)$ is continuous between Banach space we conclude that there exists a constant $\tilde\kappa_q(Q)>0$ such that 
\begin{align}||B||_{W^{1,q}}(Q)\le\tilde\kappa_q(Q)||\beta||_{W^{1-\frac{1}{q},q}(\partial Q)}.\end{align}
Since we already know that $\beta\in\Omega_{W^{1,\infty}}^1(\partial Q)$, we can get rid of the $(n-2)$-skeleton of $Q$ and reduce ourselves to perform the estimate on every face. Fix an open face $F\subset\partial Q$ of $\partial Q$ and consider the inclusion $i_F:F\rightarrow\partial Q$. Since $i_F^*\beta\in\Omega^1(F)$ and still satisfies the conditions (A.3) on $F$, by applying $d^*$ to the first equation, $d$ to the second equation and summing them up term by term, we get that
\[\Delta\beta=d^*f\in W^{-1,q}(F),\]
for every $q\in [1,+\infty)$. Now, we divide the two cases:
\begin{enumerate}
\item Notice that we can always pick a $\tilde q\in(1,+\infty)$ such that $\Delta^{-1}:W^{-1,1}(F)\rightarrow W^{1-\frac{1}{\tilde q},\tilde q}(F)$ is a continuous operator. Hence, since $d^*:L^1(F)\rightarrow W^{-1,1}(F)$ is a continuous linear operator, we conclude that
\[||\beta||_{W^{1-1/\tilde q,\tilde q}(F)}=||\Delta^{-1}(\Delta\beta)||_{W^{1-1/\tilde q,\tilde q}(F)}=||\Delta^{-1}(d^*f)||_{W^{1-1/\tilde q,\tilde q}(F)}\le\hat\kappa_1(F)||f||_{L^1(F)},\]
for some $\hat\kappa_1(F)>0$. Hence, we get the required estimate by summing up the contributes for each face of $\partial Q$ and using (A.4). 
\item As $\Delta^{-1}:W^{-1,q}(F)\rightarrow W^{1,q}(F)$ and $d^*:L^q(F)\rightarrow W^{-1,q}(F)$ are continuous linear operators, we get the estimate 
\[||\beta||_{W^{1,q}(F)}=||\Delta^{-1}(\Delta\beta)||_{W^{1,q}(F)}=||\Delta^{-1}(d^*f)||_{W^{1,q}(F)}\le\hat\kappa_q(F)||f||_{L^q(F)},\]
for some $\hat\kappa_q(Q)>0$. Hence, we get the required estimate by summing up the contributes for each face of $\partial Q$, noticing that $W^{1,q}(\partial Q)$ injects continuously into $W^{1-\frac{1}{q},q}(\partial Q)$ and using (A.4). 
\end{enumerate}
\end{proof}
\end{lem}
\begin{lem}
\label{LA.3}
Let $Q\subset\r^n$ be an open cube in $\r^n$ with centre $x_0$ and side-length $0<r<1$. Let $B\in\Omega_c^1(\r^n)$, $f\in L^{\infty}(\partial Q)$ and assume that the $1$-form $\beta:=i_{\partial Q}^*B\in\Omega_{W^{1,\infty}}^1(Q)$ is such that
\begin{align}\begin{cases}d\beta=f\\d^*\beta=0\end{cases}\mbox{ on } \partial Q.\end{align}
Then, there exists a unique solution $\tilde A$ of the following differential problem:
\begin{align}\begin{cases}\Delta A=0&\mbox{ on } Q\\i_{\partial Q}^*A=\beta &\mbox{ on } \partial Q.\end{cases}\end{align}
Moreover, $\tilde A\in\Omega_{W^{2,q}}^1(Q)\cap\Omega^1(Q)$, for every $q\in[1,+\infty)$. In particular:
\begin{enumerate}
\item $\tilde A\in\Omega_{C^{1,\alpha}}^1(\overline{Q})$ for every $\alpha\in[0,1)$;
\item $i_{\partial Q}^*\big(d^*\tilde A\big)=0$;
\item for every $q\in[1,+\infty)$ there exists a constant $K_q>0$ such that
\begin{align}\big|\big|d\tilde A\big|\big|_{L^q(Q)}\le K_qr^{1/q}\big|\big|f\big|\big|_{L^q(\partial Q)}.\end{align}
\end{enumerate}
\begin{proof}
First, we claim the existence of a weak solution $\bar A\in W_0^{1,2}\big(C_{\varepsilon,j}\big)$ for the differential system
\begin{align}\begin{cases}\Delta A=\Delta B &\mbox{ on } Q\\i_{\partial Q}^*A=0 &\mbox{ on } \partial Q.\end{cases}\end{align}
Consider the space $\Omega_{W_0^{1,2}}^1(Q)$ endowed with the scalar product defined by
\[(A_1,A_2)_{W_0^{1,2}}:=(dA_1,dA_2)_{L^2}+(d^*A_1,d^*A_2)_{L^2},\]
where
\[(A_1,A_2)_{L^2}:=\int_{Q}\star\big(A_1\wedge\star A_2\big)\, d\L^3.\]
For a proof of the fact that $(\cdot,\cdot)_{W_0^{1,2}}$ is a actually a scalar product on $\Omega_{W_0^{1,2}}^1(Q)$ we refer the reader to e.g. \cite[Lemma 2.4.10]{schwarz}.
Define the functional $E:\Omega_{W_0^{1,2}}^1(Q)\rightarrow\r$ given by 
\[E(A):=\frac{1}{2}(A,A)_{W_0^{1,2}}-(A,\Delta B)_{L^2}.\]
The existence of a global minimizer $\bar A\in\Omega_{W_0^{1,2}}^1\big(C_{\varepsilon,j}\big)$ for $E$ follows directly by applying the Lax-Milgram theorem to $E$ and the claim follows since every global minimizer for $E$ solves weakly the system (A.8).\\
Next,we notice that the $1$-form $\tilde A:=\bar A-B\in\Omega_{W^{1,2}}^1(Q)$ is a weak solution for the differential system (A.6). Moreover, since the Laplacian on $k$-forms in $\r^n$ with the flat metric coincides with the component-wise Laplacian, the uniqueness of the solution follows easily by maximum principle and we can apply Lemma \ref{LA.1} to each equation of the system (A.4) to conclude that for every $q\in[1,+\infty)$ it holds that $\tilde A\in\Omega_{W^{2,q}}^1(Q)\cap\Omega^1(Q)$ and we can find a constant $K_q^1>0$ such that
\begin{align}\big|\big|d\tilde A\big|\big|_{L^q(Q)}\le\big|\big|\nabla\tilde A\big|\big|_{L^q(Q)}\le K_q^1\big|\big|\nabla B\big|\big|_{L^q(Q)}\le K_q^1\big|\big|B\big|\big|_{W^{1,q}(Q)}.\end{align}
Thus, (1) follows directly by standard Sobolev embedding theorem. Since $A$ is $C^1(\overline{Q})$ and $i_{\partial Q}$ is a Lipschitz map, to obtain (2) we can perform the following direct computation:
\[i_{\partial Q}^*\big(d^*\tilde A\big)=d^*\big(i_{\partial Q}^*\tilde A\big)=d^*\beta=0,\]
where the last equality follows from (A.5). By applying the estimate (A.9) and Lemma \ref{LA.2} we get that there exists a constant $K_q(Q)>0$ such that
\[\big|\big|d\tilde A\big|\big|_{L^q(Q)}\le K_q(Q)\big|\big|f\big|\big|_{L^q(\partial Q)}.\]
Eventually, (3) follows easily by setting $K_q:=K_q(C)$ and applying a scaling argument similar to the one in the proof of \ref{LA.1} to the previous estimate.
\end{proof}
\end{lem}
\begin{dfn}
\label{DA.1}
The $1$-form $\tilde A$ whose existence has been proved by the previous Lemma \ref{LA.3} is called the \textbf{harmonic extension} of the boundary datum $\beta$ to the open cube $Q$. 
\end{dfn}
\section{A useful characterization of divergence free vector fields in $L_{loc}^1(\Omega)$}
\begin{lem}
\label{LB.1}
Let $\Omega\subset\r^n$ be any open subset of $\r^n$ and $X\in L_{loc}^1(\Omega,\r^3)$ be any vector field over $\Omega$. Then the following are equivalent:
\begin{enumerate}
\item $\Div(X)=0$ in $\mathcal{D}'(\Omega)$;
\item for every Lipschitz and compactly supported function $\varphi:\Omega\rightarrow\r$ it holds that
\[\int_{\Omega}X(x)\cdot\nabla\varphi(x)\,dx=0.\]
\end{enumerate}
\begin{proof}
As the implication $2\Rightarrow 1$ is trivial, it suffices to show that $1\Rightarrow 2$.\\ Since $\varphi:\Omega\rightarrow\r$ is Lipschitz, we know that $\nabla\varphi\in\ L^{\infty}(\Omega,\r^n)$. For every $\varepsilon>0$, define the set
\[\Omega_{\varepsilon}:=\{x\in\Omega \mbox{ s.t. } \dist(x,\partial\Omega)>\varepsilon\}.\]
Since $\supp(\varphi)$ is closed in $\Omega$, we can pick some $\varepsilon_0>0$ such that $\Omega_{\varepsilon_0}\supset\supp(\varphi)+\overline{B_{\varepsilon_0}(0)}$. Consider a regularizing kernel $\{\rho_{\varepsilon}\}_{0<\varepsilon<\varepsilon_0}\subset C_c^{\infty}(\r^n)$ such that $\supp(\rho_{\varepsilon})\subset B_{\varepsilon}(0)$, for every $0<\varepsilon<\varepsilon_0 $. Define the family $\{\varphi_{\varepsilon}\}_{0<\varepsilon<\varepsilon_0}\subset C_c^{\infty}(\Omega_{\varepsilon})$ by
\[\varphi_{\varepsilon}(x):=(\varphi\star\rho_{\epsilon})(x)=\int_{B_{\varepsilon}(x)}\varphi(y)\rho_{\varepsilon}(x-y)\, dy \quad \mbox{ for every } x\in\Omega_{\varepsilon}.\] 
We observe that the family $\left\{\nabla\varphi_{\varepsilon}|_{\Omega_{\varepsilon_0}}\right\}_{0<\varepsilon<\varepsilon_0}$ is bounded in $L^{\infty}(\Omega_{\varepsilon_0})$. Indeed,
\[\left|\nabla\varphi_{\varepsilon}|_{\Omega_{\varepsilon_0}}(x)\right|\le\int_{B_{\varepsilon}(x)}|\nabla\varphi(y)\rho_{\varepsilon}(x-y)|\, dy\le||\nabla\varphi||_{L^{\infty}(\Omega)}, \quad \mbox{ for every } x\in\Omega_{\varepsilon_0}.\]
Thus, as bounded subsets of $L^{\infty}(\Omega_{\varepsilon_0})$ are $\text{weak}^*$ compact, we can find a sequence $\left\{\nabla\varphi_{n_k}|_{\Omega_{\varepsilon_0}}\right\}_{k\in\n}$ such that $\nabla\varphi_{n_k}|_{\Omega_{\varepsilon_0}}\stackrel{\ast}{\rightharpoonup} V$ for some $V\in L^{\infty}(\Omega_{\varepsilon_0})$. As $\varphi|_{\Omega_{\varepsilon_0}}\in L^1(\Omega_{\varepsilon_0})$, it holds that $\varphi_{n_k}|_{\Omega_{\varepsilon_0}}\rightarrow \varphi|_{\Omega_{\varepsilon_0}}$ strongly in $L^1(\Omega_{\varepsilon_0})$ and thus also in $\mathcal{D}'(\Omega_{\varepsilon_0})$. Since the weak gradient operator is continuous with the respect to the topology on $\mathcal{D}'(\Omega_{\varepsilon_0})$, we conclude that $\nabla\varphi_{n_k}|_{\Omega_{\varepsilon_0}}\rightarrow\nabla\varphi|_{\Omega_{\varepsilon_0}}$ in $\mathcal{D}'(\Omega_{\varepsilon_0})$ and thus $V=\nabla\varphi|_{\Omega_{\varepsilon_0}}$ in $L^{\infty}(\Omega_{\varepsilon_0})$. Hence, we deduce that
\[\int_{\Omega_{\varepsilon_0}}\tilde X(x)\cdot\nabla\varphi_{n_k}(x)\, dx\rightarrow\int_{\Omega_{\varepsilon_0}}\tilde X(x)\cdot\nabla\varphi(x)\, dx,\]
for every $X\in L^1(\Omega_{\varepsilon_0})$. By picking $\tilde X:=\rchi_{\supp(\varphi)+\overline{B_{\varepsilon_0}(0)}}X|_{\Omega_{\varepsilon_0}}$ in the previous equality we get that
\[0=\int_{\Omega}X(x)\cdot\nabla\varphi_{n_k}(x)\, dx\rightarrow\int_{\Omega}X(x)\cdot\nabla\varphi(x)\, dx\]
and the statement follows.
\end{proof}
\end{lem}
\begin{lem}
\label{LB.2}
Let $\Omega\subset\r^n$ be any open subset of $\r^n$ and $X\in L_{loc}^1(\Omega,\r^3)$ be any vector field over $\Omega$. Then the following are equivalent:
\begin{enumerate}
\item $\Div(X)=0$ in $\mathcal{D}'(\Omega)$;
\item for every $x_0\in\Omega$ and for a.e. $\displaystyle{r\in\left(0,\frac{1}{\sqrt{n}}\dist(x_0,\partial\Omega)\right)}$ it holds that
\[\int_{\partial C_r(x_0)}\big(X\cdot\nu_{\partial C_r(x_0)}\big)\,d\H^{n-1}=0,\]
where $\nu_{\partial C_r(x_0)}:\partial C_r(x_0)\rightarrow S^2$ is the unit outward pointing normal vector to $\partial C_r(x_0)$.
\end{enumerate} 
\begin{proof}
We begin by proving that $1\Rightarrow 2$.
Fix any $x_0\in\Omega$ and define $\displaystyle{d:=\frac{1}{\sqrt{n}}\dist(x_0,\partial\Omega)}$. Look at the function $g:(0,d)\rightarrow\r$ defined as
\[g(s):=\int_{\partial C_s(x_0)}\big(X\cdot\nu_{\partial C_s(x_0)}\big)\, d\H^{n-1}, \quad \mbox{ for every } s\in(0,d). \]
Notice that, since $X\in L_{loc}^1(\Omega)$, the function $h:(0,d)\times\Omega\rightarrow\r$ given by
\[h(s,x):=\left(X(x)\cdot\nu_{\partial C_s(x_0)}(x)\right)\rchi_{\partial C_s(x_0)}(x), \quad \mbox{ for eveery } (s,x)\in(0,d)\times\Omega\] 
belongs to $L^1\left((0,d)\times\Omega, \L^1|_{(0,d)}\times\H^{n-1}|_{\Omega}\right)$. Thus, by Fubini's theorem, $g$ is a well-defined function in $L^1\left((0,d),\L^1|_{(0,d)}\right)$. Then, by the Lebesgue points theorem, a.e. $s\in(0,d)$ is a Lebesgue point of $g$. Pick $r\in(0,d)$ to be any Lebesgue point of $g$.
For every sufficiently small $0<\varepsilon<\varepsilon_0<r$, where $\varepsilon_0>0$ is chosen such that $C_{r+\varepsilon_0}(x_0)\subset\Omega$, define the Lipschitz function $f_{\varepsilon}:[0,+\infty)\rightarrow[0,1]$ by
\[f_{\varepsilon}(t):=\begin{cases}1 \ &\mbox{ if } t\in [0,r-\varepsilon),\\\displaystyle{\frac{r+\varepsilon}{2\varepsilon}-\frac{t}{2\varepsilon}} \ &\mbox{ if } t\in [r-\varepsilon,r+\varepsilon),\\0 \ &\mbox{ if } t\in [r+\varepsilon,+\infty).\end{cases}\]
and use it build the compactly supported Lipschitz function $\varphi_{\varepsilon}:\Omega\rightarrow\r$ given by
\[\varphi_{\varepsilon}(x):=f_{\varepsilon}\big(||x-x_0||_{*}\big), \quad \mbox{ for every } x\in\Omega,\]
where
\[||x||_{*}:=\sup_{j=1,...,n}|x_j|, \qquad \mbox{ for every } x=(x_1,...,x_n)\in\r^n.\]
Now, by \cite[Lemma 7.6]{maggi}, since $f_{\varepsilon}$ is piecewise affine, differentiable on $[0,+\infty)\smallsetminus F:=\{r-\varepsilon,r+\varepsilon\}$ and $u:=||\cdot - \ x_0||_{*}\in W_{loc}^{1,\infty}(\Omega)\subset W_{loc}^{1,1}(\Omega)$, then $\nabla u=0$ a.e. on $u^{-1}(F)=\partial C_{r-\varepsilon}(x_0)\cup \partial C_{r+\varepsilon}(x_0)$ and $\varphi_{\varepsilon}=f_{\varepsilon}\circ u\in W_{loc}^{1,1}(\Omega)$ with
\[\nabla\varphi_{\varepsilon}(x):=f_{\varepsilon}'\big(||x-x_0||_{*}\big)\nabla u(x)=-\frac{1}{2\varepsilon}\rchi_{C_{r+\varepsilon}(x_0)\smallsetminus C_{r-\varepsilon}(x_0)}(x)\nu_{\partial C_{||x-x_0||_{*}}(x_0)}(x), \quad \mbox{ for a.e. } x\in\Omega.\]
By the previous Lemma \ref{LB.1}, it holds that
\[\int_{\Omega}X(x)\cdot\nabla\varphi_{\varepsilon}(x)\, dx=0,\]
for every $0<\varepsilon<\varepsilon_0$. Moreover, since $r$ is a Lebesgue point of $g$, we know that
\[\int_{\Omega}X(x)\cdot\nabla\varphi_{\varepsilon}(x)\, dx=-\frac{1}{2\varepsilon}\int_{r-\varepsilon}^{r+\varepsilon}\int_{\partial C_s(x_0)}X(y)\cdot\nu_{\nu_{\partial C_s(x_0)}}(y)\, d\H^{n-1}(y)ds=-\frac{1}{2\varepsilon}\int_{r-\varepsilon}^{r+\varepsilon}g(s)\, ds\rightarrow -g(r)\]
as $\varepsilon\rightarrow 0^{+}$. Thus, by uniqueness of the limit,
\[g(r)=\int_{\partial C_r(x_0)}\big(X\cdot\nu_{\partial C_r(x_0)}\big)\, d\H^{n-1}=0\]
and $2$ follows.\\
We move now to show that $2\Rightarrow 1$. We claim that
\[\int_{\Omega}X(x)\cdot\nabla\varphi(x)\,dx=0, \quad \mbox{ for every } ||\cdot||_{*}\mbox{-radial function } \varphi\in \Lip_c(\Omega).\]
Indeed, by definition, $\varphi\in\Lip_c(\Omega)$ is $||\cdot||_{*}|$-radial with respect to some point $x_0\in\Omega$ if there exists $f\in\Lip_c\left([0,+\infty)\right)$ with $\supp(f)\subset[0,r_0)$ for some $0<r_0<\dist(x_0,\partial\Omega)$ such that $\varphi=f\left(||\cdot- \ x_0||_{*}\right)$.  Then, it holds that
\[\begin{aligned}\int_{\Omega}X(x)\cdot\nabla\varphi(x)\, dx&=\int_{\Omega}f'\left(||x-x_0||\right)X(x)\cdot\nu_{\partial C_{||x-x_0||_{*}}(x_0)}(x)\, dx=\\
&=\int_{0}^{r_0}f'\left(r\right)\int_{\partial Cr(x_0)}X(x)\cdot \nu_{\partial C_r(x_0)}(x)\, d\H^{n-1}(x)=0.
\end{aligned}\]
For every $\varepsilon>0$, define the set
\[\Omega_{\varepsilon}:=\{x\in\Omega \mbox{ s.t. } \dist(x,\partial\Omega)>\sqrt{n}\varepsilon\}\]
and consider a $||\cdot||_{*}$-radial and positive regularizing kernel $\{\rho_{\varepsilon}\}_{\varepsilon>0}\subset \Lip_c(\r^n)$ such that,  for every $\varepsilon>0$, $\supp(\rho_{\varepsilon})\subset C_{\varepsilon}(0)$. Define the $\varepsilon$-regularization $\{X_{\varepsilon}\in\Lip(\Omega_{\varepsilon},\r^n)\}_{\varepsilon>0}$ of $X$ by
\[X_{\varepsilon}(x):=(X\star\rho_{\epsilon})(x)=\int_{C_{\varepsilon}(x)}X(y)\rho_{\varepsilon}(x-y)\, dy \quad \mbox{ for every } x\in\Omega_{\varepsilon}.\] 
We claim that $X_{\varepsilon}$ is divergence free, for every $\varepsilon>0$. Indeed, for every $x\in\Omega_{\varepsilon}$, it holds that
\[\Div(X_{\varepsilon})(x)=\int_{C_{\varepsilon}(x)}X(y)\cdot\nabla\rho_{\varepsilon}(x-y)\, dy=\int_{\Omega}X(y)\cdot\nabla\rho_{\varepsilon}(x-y)\, dy=0,\]
by exploiting the previous claim, since $\rho_{\varepsilon}(x-\cdot)|_{\Omega}\in\Lip_c(\Omega)$ is a radial function with respect to the point $x\in\Omega_{\varepsilon}\subset\Omega$. Fix any $\varphi\in C_c^{\infty}(\Omega)$ and find $\varepsilon_0>0$ such that $\supp(\varphi)\subset\Omega_{\varepsilon_0}$. As $X|_{\Omega_{\varepsilon_0}}\in L^1_{loc}(\Omega_{\varepsilon_0})$, then $X_{\varepsilon}|_{\Omega_{\varepsilon_0}}\rightarrow X|_{\Omega_{\varepsilon_0}}$  in $\mathcal{D}'(\Omega_{\varepsilon_0})$. Thus, for every $0<\varepsilon<\varepsilon_0$, it holds that
\[0=\left<\Div(X_{\varepsilon}),\varphi\right>=-\int_{\Omega_{\varepsilon_0}}X_{\varepsilon}(x)\cdot\nabla\varphi(x)\, dx\rightarrow -\int_{\Omega_{\varepsilon_0}}X(x)\cdot\nabla\varphi(x)\, dx=-\int_{\Omega}X(x)\cdot\nabla\varphi(x)\, dx\]
and this concludes the proof.
\end{proof}
\end{lem}
\printbibliography
\end{document}